\documentclass[11pt,a4paper]{article}
\usepackage{graphicx}
\usepackage{latexsym,color}
\usepackage{amssymb,amsmath}
\newcommand{\norm}[2]{
\left\| #2 \right\|_{#1}
}
\newcommand{\Hil}[0]{
\mathcal{H} 
}
\newcommand{\BL}[0]{
{\mathcal B}
}

\newcommand{\RR}[0]{\mathbb{R} }
\newcommand{\BB}[0]{{\mathfrak B}}

\newcommand{\CC}[0]{\mathbb{C}}

\newcommand{\MM}[0]{\mathbf {M}}

\newtheorem{theorem}{Theorem}[section]
\newtheorem{def.}[theorem]{Definition}
\newtheorem{prop.}[theorem]{Proposition}
\newtheorem{lem.}[theorem]{Lemma}
\newtheorem{cor.}[theorem]{Corollary}
\newtheorem{conj.}[theorem]{Conjecture}
\newtheorem{Bsp.}{Example}[section]
\newtheorem{rem.}{Remarks}[theorem]

\renewcommand{\ge}{\geqslant}
\renewcommand{\le}{\leqslant}

 \newcommand{\noi}{\noindent}
\newcommand{\dis}{\displaystyle}

 \newcommand{\bea}{\begin{eqnarray}}
\newcommand{\ena}{\end{eqnarray}}

\newcommand{\beano}{\begin{eqnarray*}}
\newcommand{\enano}{\end{eqnarray*}}
 
\newcommand{\bei}{\begin{itemize}}
\newcommand{\eni}{\end{itemize}}

\newcommand{\bedefi}{\begin{def.} \rm }
\newcommand{\eofproof}{\hfill $\square$\bigskip} 
\newenvironment{proof}{\noindent \bf Proof: \rm}{\eofproof

\vspace{5mm}}

\begin{document}
 
\markboth{Balazs, Antoine and Grybos}{Weighted and Controlled Frames}

\title{Weighted and Controlled Frames}
\author{Peter Balazs, \\ {\small Acoustics Research Institute, Austrian Academy of Sciences,} \\ {\small Wohllebengasse 12-14, 
A-1040 Vienna, Austria}\\{\small Peter.Balazs@oeaw.ac.at}\\
Jean-Pierre Antoine, \\ {\small Institut de Physique Th\'eorique, Universit\'e Catholique de Louvain,} \\{\small chemin du Cyclotron 2, 
B-1348 Louvain-la-Neuve, Belgium}\\{\small antoine@fyma.ucl.ac.be}
\\
Anna Grybo\'s, \\ {\small NuHAG, Faculty of Mathematics, University of Vienna,} \\ {\small Nordbergstrasse 15, A-1090 Austria;} \\
 {\small Faculty of Applied Mathematics, AGH University of Science and Technology,} \\ {\small Al. Mickiewicza 30, 30-059 Krakow, Poland}\\
 {\small anna.grybos@univie.ac.at}
}

\maketitle

 \begin{abstract} Weighted and controlled frames have been introduced recently to improve the numerical efficiency of 
 iterative algorithms for inverting the frame operator. In this paper we develop systematically these notions, 
  including their mutual relationship. We will show that controlled frames are equivalent to standard frames and so this concept gives a 
  generalized way to check the frame condition, while offering a numerical advantage  in the sense of preconditioning. 
Next, we investigate weighted frames, in particular their relation to controlled frames. 
  We consider the special case of semi-normalized weights, where the concepts of weighted frames and standard frames are interchangeable. 
 We also make the connection with frame multipliers. 
 Finally we analyze weighted frames numerically. First we investigate three possibilities for finding weights in order to tighten a given frame, i.e., 
decrease the frame bound ratio. Then we examine Gabor frames and
 how well the canonical dual of a weighted frame is approximated by the inversely weighted dual frame.  
\end{abstract}



\section{Introduction}

In practice, the frame bound ratio of a given nontight frame can often be reduced by weighting the elements. 
The so-called weighted frames, i.e.,  frames $(\psi_n)$ with complex weights $(\omega_n)$ such that the sequence $(\omega_k   \psi_k)$ is again a frame, 
were introduced in Ref.\cite{bogdvan1} to get a numerically more efficient approximation algorithm for spherical wavelets. 
 By decreasing the ratio of the frame bounds, weighting improves the numerical efficiency of iterative algorithms like the `{\em frame algorithm}'\cite{ole1}
 for the inversion of the frame operator. The same paper\cite{bogdvan1} introduced and used controlled frames, that is, a frame $(\psi_n)$ and an operator $C$ 
 such that the combination of $C$ with the frame operator $L$ is positive and invertible. 
Since these concepts were used there just as a tool for spherical wavelets, they were not discussed in full detail.  

In this paper, we will develop the related theory and derive some properties used in Ref.\cite{bogdvan1} without proof, as well as give the results of numerical experiments. 
 Section \ref{sec:prelimn0} contains some preliminary results.  
In Section \ref{sec:contrfram1} we will show that controlled frames are equivalent to standard frames and so this concept gives a generalized way to check the frame condition. 
In Section \ref{sec:weightfram0} we investigate weighted frames.
We will put some emphasis on the mutual relationship between the two concepts, showing in particular that weighted frames cannot always be considered as controlled frames.
We also investigate how these concepts can improve the efficiency of iterative algorithms for inverting the frame operator. As a special case, 
we consider weights bounded and bounded away from zero, for which the concepts of frames and weighted frames are interchangeable again. 
The connection to frame multipliers will be addressed briefly.

In the last part we will investigate the concept of weighted frames in numerical experiments. In Section \ref{sec:weightnum0} we analyze three different
choices for weights with the aim of making frames tighter, i.e.,  reducing the quotient of the frame bounds. 
We give the results of some numerical experiments, showing that these weights  very often improve the condition number of the frame operator matrix.
We see that redundancy is an important parameter for the optimality of these weights. In Section \ref{sec:gabornum0} we examine the computational behavior
 of weighted Gabor frames. In particular we investigate how well the canonical dual weighted frame is approximated by the inversely weighted dual frame.
 We see that the error depends linearly on %the range of the weights and 
the amount of weighted elements and the redundancy. 

\vspace*{-2mm}
\section{Preliminaries}\label{sec:prelimn0}

In this section, we collect the basic notation and some preliminary results. 
Throughout the paper, $\Hil$ is a separable Hilbert space, with inner product  $\left< ., . \right>$,  linear in the first coordinate, and norm $\norm{}{\cdot}$. 
We denote by $I$ the identity operator on  $\Hil$. Let $\BL(\Hil_1,\Hil_2)$ be the set of all bounded linear operators from $\Hil_1$ to $\Hil_2$.  
This set is a Banach space for the   operator norm   $\norm{}{A} = \dis\sup _{\norm{\Hil_1}{x} \le 1}   \norm{\Hil_2}{A x}  $.
 The  {adjoint} of the operator $A$ is denoted by $A^*$ and the spectrum of $A$ by $\sigma(A)$.  We define   $GL(\Hil_1,\Hil_2)$ as  
the set of all bounded linear operators with a bounded inverse, and similarly for $GL(\Hil)$.  
 Our standard reference for  Hilbert space and operator theory is Ref.\cite{conw1}.

\vspace*{-4mm}
\subsection{Frames} 

We collect here some known facts about frames, in a form suitable for us.
For more details on this topic,  see for instance Refs.\cite{Casaz1} or  \cite{ole1}. 
\bedefi \label{sec:framprop1} 
A sequence $\Psi = (\psi_n, n \in \Gamma)$ is called a 
\em frame \em for the Hilbert space $\Hil$, if there exist constants ${\rm\sf m} > 0 $ and ${\rm\sf M}<\infty$  such that \\
\centerline{${\rm\sf m}  \norm{}{f}^2 \le \sum \limits_{n\in \Gamma} \left| \left< f, \psi_n \right> \right|^2 \le {\rm\sf M}  
\norm{}{f}^2 ,  \forall \, f \in \Hil.$ }
${\rm\sf m}$ is a \em lower\em , ${\rm\sf M}$ an \em upper frame bound\em .
If the bounds can be chosen such that ${\rm\sf m}={\rm\sf M}$, the frame is called \em tight\em .  
\end{def.}
The optimal bounds ${\rm\sf m}_{opt}, {\rm\sf M}_{opt}$ are the largest ${\rm\sf m}$ and smallest ${\rm\sf M}$ that fulfill the corresponding inequality.

\bedefi 
Given a frame $\Psi = (\psi_n, n \in \Gamma)$,  
$L_\Psi%{( \psi_n )}
 : \Hil  \rightarrow \Hil $ denotes the \em (associated) frame  operator \em
$L_{\Psi} ( f  ) = \sum \limits_n  \left< f , \psi_n \right>  \psi_n $.
%and $W_{(\psi_n)}: l^2  \rightarrow \Hil$  
%the (associated) synthesis operator $W_{(\psi_n)} (c_n) = \sum \limits_n c_n \psi_n$.
\end{def.}
If there is no risk of confusion, we will omit the index and write $L$ instead of $L_{\Psi}$. 
For any frame,  $L$   is a  positive   invertible operator on all of $\Hil$, satisfying the inequalities ${\rm\sf m}\, I  \le L \le {\rm\sf M}\; I $ 
and ${\rm\sf M}^{-1} I  \le L^{-1} \le {\rm\sf m}^{-1} I $. % As a matter of fact, frames can be classified by this property. 
Furthermore,

\begin{theorem} 
Let $\Psi = \left( \psi_n \right)$ be a frame for $\Hil$ with  bounds ${\rm\sf m}, {\rm\sf M} > 0$. 
Then $\tilde \Psi = \left( \widetilde{\psi}_n \right) = \left( L^{-1} \psi_n \right)$ is a frame with bounds ${\rm\sf M}^{-1}, {\rm\sf m}^{-1} > 0$, 
the so-called \em canonical dual frame\em . %\index{frames!canonical dual}.
 Every $f \in \Hil$ has  expansions
$ f = \sum \limits_{n \in \Gamma} \left< f, \widetilde \psi_n \right> \psi_n $
and 
$ f = \sum \limits_{n \in \Gamma} \left< f, \psi_n \right> \widetilde \psi_n $,
and both sums converge unconditionally in $\Hil$.
\end{theorem}

We will use so-called {\em frame multipliers}.\cite{xxlmult1} These are operators defined by
$$ 
\MM_{{\bf m}, \Psi, \Phi} f = \sum \limits_k m_k \left< f , \psi_k \right> \phi_k $$ 
for the frames $\Psi=(\psi_n)$ and $\Phi=(\phi_n)$ and the weight sequence ${\bf m}=(m_k)$. We shorten the notation by setting $\MM_{{\bf m}, \Psi} = \MM_{{\bf m}, \Psi, \Psi}$. 

%\subsubsection{Gabor frames} \label{sec:gaborframes}
Among all frames, a privileged role is played by Gabor\cite{Groech1} and wavelet frames.\cite{daubech1} 
For future use, let us repeat the definition of the former.

Given $a,b>0$, a \emph{Gabor frame} over the   regular lattice $\Lambda = a \mathbb{Z}^d \times b \mathbb{Z}^d $ is a family 
$$ 
G = (g_{k,l})_{k,l}:= \{ M_{lb}T_{ka}g, \;l, k \in \mathbb{Z}^{d} \}
$$ 
that fulfills the frame condition and whose elements are translated and modulated versions of a given window function $g\in L^2(\RR^d)$. Here the operations 
of translation  $T_x$ and modulation $M_\xi$ are defined by: 
$$ 
T_x f(t)=f(t-x) \ \mathrm{and} \ M_\xi f(t)=e^{2 \pi i\xi t} f(t),\; t,x,\xi \in  \RR^d. 
$$
We will denote as $\lambda =(\tau, \omega) \in \Lambda$ the time-frequency shift, defined as
$$ \pi (\lambda) g = M_\omega T_\tau g. $$ 
We shall discuss some concrete examples of (weighted, discrete) Gabor
 frames in Section \ref{sec:gabornum0}.

\vspace*{-2mm}
\subsection{The bounded and boundedly invertible positive operators $GL^{(+)}(\Hil)$}

A bounded operator $T$ is called \em positive \em (respectively \em non-negative\em), if \mbox{$\left< T f , f\right> > 0$} 
for all $f \not=0$ (respectively $\left< Tf, f \right> \ge 0$ for all $f$). Every non-negative operator is clearly self-adjoint. 
If $A \in \BL(\Hil)$ is non-negative, then there exists  a unique non-negative operator $B$ such that $B^2 = A$.\cite{gohbgol1}
 Furthermore $B$ commutes 
with every operator that commutes with $A$. This will be denoted by $B = A^{1/2}$.
Let $GL^{(+)}(\Hil)$ be the set of positive operators in $GL(\Hil)$.
 
The following result is needed in the sequel, but straightforward to prove:
\begin{prop.} \label{sec:posinvertop1} 
Let $T : \Hil \rightarrow \Hil$ be a linear operator. Then the following conditions are equivalent:
\begin{enumerate}
\item There exist $m > 0$ and $M < \infty$, such that ${\rm\sf m}\,  I \le T \le {\rm\sf M}\, I$;
\item $T$ is positive and there exist $m > 0$ and $M < \infty$, such that ${\rm\sf m} \norm{}{f}^2 \le \norm{}{T^{1/2}f}^2 \le {\rm\sf M} \norm{}{f}^2$;    
 \item $T$ is positive and $T^{1/2} \in GL(\Hil)$;
\item There exists a self-adjoint operator $A \in GL(\Hil)$,  such that $A^2=T$;
\item $T \in GL^{(+)}(\Hil)$;
\item There exist constants ${\rm\sf m} > 0$ and ${\rm\sf M} < \infty$ and an operator $C \in GL^{(+)}(\Hil)$ such that ${\rm\sf m}' C\le T \le {\rm\sf M}' C$; 
\item For every $C \in GL^{(+)}(\Hil)$, there exist constants ${\rm\sf m} > 0$ and ${\rm\sf M} < \infty$ such that ${\rm\sf m}' C \le T \le {\rm\sf M}' C$. 
\end{enumerate}
\end{prop.}

\bedefi Given $T \in GL^{(+)}(\Hil)$, any two constants ${\rm\sf m}_{T},{\rm\sf M}_{T}$ such that 
$${\rm\sf m}_{T}\, I  \le T \le {\rm\sf M}_{T}\,  I $$
 are called \em lower \em and \em upper bound of $T$\em, respectively. If ${\rm\sf m}_T$ is maximal, 
resp. if ${\rm\sf M}_T$ is minimal, we call them the \em optimal \em bounds and we denote them by ${\rm\sf m}_T^{(opt)}, {\rm\sf M}_T^{(opt)}$.
\end{def.}
The upper and lower bounds are clearly not unique. 

The following results are easily proved using Proposition \ref{sec:posinvertop1}:
\begin{cor.}  
Let $T \in GL^{(+)}(\Hil)$. Then 
\begin{enumerate}
\vspace*{-2mm}
\item $\left\| T \right\|_{} = {\rm\sf M}_T^{(opt)}$.
\item $\sigma(T) \subseteq \left[ {\rm\sf m}_T , {\rm\sf M}_T \right]$, for any lower, resp. upper, bounds. 
\end{enumerate}
\end{cor.} 

\begin{cor.} \label{sec:inversposbounds1} For $T \in GL^{(+)}(\Hil)$, the numbers ${\rm\sf m}_{T^{-1}} = {\rm\sf M}_T^{-1}$ and 
${\rm\sf M}_{T^{-1}} = {\rm\sf m}_T^{-1}$ are bounds for $T^{-1}$. In particular $\left\| T^{-1} \right\|_{} =1/ {{\rm\sf m}_T^{(opt)}}.$
\end{cor.}

\begin{cor.} \label{sec:commutposoperbound1} Let $S,T \in GL^{(+)}(\Hil)$ be commuting operators. Then 
$T$ admits as lower and upper bounds $(\frac{{\rm\sf m}_{TS}}{{\rm\sf M}_S}, \frac{{\rm\sf M}_{TS}}{{\rm\sf m}_S})$ and $S \, T$ admits 
$({\rm\sf m}_S {\rm\sf m}_T, {\rm\sf M}_S {\rm\sf M}_T)$.
\end{cor.}

\subsection{Numerical issues}  \label{sec:numiss0} 

A well-known algorithm to find the inverse of an operator is the Neumann algorithm, which is based on the following property:
\begin{prop.} \label{sec:neumanalgor1} 
Given two Banach spaces  $\BB_1,\BB_2$, if  $U:\BB_1 \rightarrow \BB_2$ is bounded and $\norm{}{I - U} < 1$, then $U$ is invertible and 
$ U^{-1} = \sum \limits_{n=0}^{\infty} \left( I - U \right)^k .$
Furthermore $\norm{}{U^{-1}} \le  (1 - \norm{}{I-U})^{-1} $.
\end{prop.}

A way to improve the numerical efficiency of an iterative algorithm for solving a linear system of equations is preconditioning.\cite{xxlfei1,trebau1}
Instead of solving the linear system of equations 
$A x = b$, 
one solves the system 
$ P A x = P b$
  for a properly chosen \emph{preconditioning} matrix $P$.

A `clustered spectrum' yields a fast convergence  as well as guarantees a small condition number,\cite{luen} 
 $\kappa\left(A\right)=\left\| A^{-1}\right\|_{} \left\| A \right\|_{}$, since
$\kappa(A) =  {\sigma_n}/{\sigma_1}$ where $\sigma_n$ and $\sigma_1$ are
the largest and smallest singular values, respectively. 
\\

Given an operator $T \in GL^{(+)}(\Hil)$ its condition number is given by 
$$ \kappa(T) = \left\| { T }^{-1} \right\|\,  \left\| T \right\| = \frac{{\rm\sf M}_{ T}^{(opt)}}{{\rm\sf m}_{T}^{(opt)}}.
 $$
So we can use preconditioning by looking for an operator $C$  
such that 
$$
 \kappa(C \, T) = \dis\frac{{\rm\sf M}_{C \, T}}{{\rm\sf m}_{C \, T}} < \dis\frac{{\rm\sf M}_{T}}{{\rm\sf m}_{T}} = \kappa(T).
 $$
 
\section{Controlled frames} \label{sec:contrfram1}

\bedefi
\label{contrframe}Let $C \in GL(\Hil)$. A \em frame controlled by the operator $C$ \em or \emph{$C$-controlled frame} %(?????) 
is a family of vectors $\Psi = \left( \psi_n \in \Hil : n \in \Gamma \right)$, 
such that there exist two constants ${\rm\sf m}_{CL} > 0 $ and ${\rm\sf M}_{CL} <\infty$ satisfying
\begin{equation} \label{sec:contrframineq1}
{\rm\sf m}_{CL} \norm{}{f}^2 \le \sum \limits_n \left< \psi_n , f \right> \left< f , C \psi_n \right> \le {\rm\sf M}_{CL} \norm{}{f}^2 , 
\quad \mbox{for all}\;  f \in \Hil.
\end{equation}
We call 
$$L_C f = \sum \limits_{n \in \Gamma} \left< \psi_n , f \right> C \psi_n $$
the \em controlled frame operator\em .
\end{def.}
The definition above is clearly equivalent to $CL\in GL^{(+)}(\Hil)$, so  the notation is coherent with the one in the previous section.
\begin{prop.} \label{sec:controlframclass1} Let $\Psi$ be a $C$-controlled frame in $\Hil$ for $C \in GL(\Hil)$. Then $\Psi$ is a classical frame.
 Furthermore $C \, L = L \, C^*$ and so
$$
 \sum \limits_{n \in \Gamma} \left< \psi_n , f \right> C \psi_n = \sum \limits_{n \in \Gamma} \left< C \psi_n , f \right>\psi_n .
$$
\end{prop.}
\noi{\bf Proof.} 
Let $\Psi$ be a controlled frame. Then using the definition and Proposition \ref{sec:posinvertop1}, we know that $L_C \in GL(\Hil)$. 
Let $ \widetilde{L} = C^{-1} L_C$. Clearly $\widetilde{L} \in GL(\Hil)$ and
$$
 \widetilde{L} f = C^{-1} \left( \sum \limits_{n \in \Gamma} \left< \psi_n , f \right> C \psi_n \right)  =
  \sum \limits_{n \in \Gamma} \left< \psi_n , f \right> \psi_n = L f
$$
Therefore $L$ is everywhere defined and $L \in GL(\Hil)$. Thus  $\Psi$ is a frame \cite{ole1}. 

By definition $L_C$ is positive, therefore  self-adjoint. So 
$L_C = C \, L = L_C^* = L^* \, C^* = L \, C^*  
$ .\vspace*{-1.4cm} 
\eofproof
\vspace{1.4cm}

\noi Since every controlled frame is a (classical) frame, \eqref{sec:contrframineq1} yields a criterion
 to check if a given sequence constitutes a frame. Furthermore it becomes obvious from the last result that the role of $C \psi_n$ and $\psi_n$ could have been switched in the definition of controlled frames.

It is difficult to see   in  full generality which conditions are needed  for a frame and an operator to form a controlled frame.
 But if $C$ is self-adjoint, we can give necessary   and sufficient conditions:
\begin{prop.} \label{sec:controlframclass2}
 Let $C \in GL(\Hil)$ be self-adjoint. The family $\Psi$ is a frame 
  controlled by $C$  
if and only if it is a (classical) frame for $\Hil$, and $C$ is positive and commutes with the frame operator $L$.
\end{prop.}
\noi{\bf Proof.}  Suppose $C$ and $\Psi$ form a controlled frame. Then from Proposition \ref{sec:controlframclass1}, it is clear that $\Psi$ is a frame and that $L$ and $C$ commute. 
Therefore $C = L_C \, L^{-1}$ is also positive.

For the converse implication, we note that, if $\Psi$ is a frame, then $L \in GL^{(+)}(\Hil)$. Therefore $C  L = L_C \in GL^{(+)}(\Hil)$ and so 
$L_C$ is positive. By Proposition \ref{sec:posinvertop1},  Eq. \eqref{contrframe} is satisfied. 
\eofproof

Using Propositions \ref{sec:controlframclass2} and \ref{sec:posinvertop1} the following result is easy to show:
\begin{cor.}  \quad
\bei
\item[(1)] Let $C$ be an invertible, self-adjoint operator and $L$ be the frame operator of $\Psi$. Then, 
$ {\rm\sf m}\, I \le C   L \le{\rm\sf M}\,  I  \; \mbox{ implies } \;{\rm\sf m}\, C^{-1} \le L \le  {\rm\sf M}\, C^{-1}.$
\item[(2)]  Let $C \in GL^{(+)}(\Hil)$. Then,
$ {\rm\sf m}\, C^{-1} \le L \le  {\rm\sf M}\, C^{-1}  \; \mbox{ implies } \; {\rm\sf m}\, I \le C   L \le{\rm\sf M}\,  I  $
and  $\Psi$ is a frame.
\eni
\end{cor.} 
%\noi{\bf Proof.}  
%(1) The inequalities on the left mean that $\Psi$ is a frame controlled by $C$, thus $C$ is positive and commutes with $L$. Therefore,
%\begin{align*}
%{\rm\sf m}\, I \le C   L \le{\rm\sf M}\,  I  & \Longleftrightarrow  \norm{}{\sqrt{{\rm\sf m}}x} 
%\le  \norm{}{{\left(CL\right)}^{1/2}x} \le \norm{}{\sqrt{{\rm\sf M}}x}, \; \forall\, x \in \Hil \\
%& \Longleftrightarrow \norm{}{\sqrt{\rm\sf m}x} \le  \norm{}{L^{1/2} C^{1/2}x} \le \norm{}{\sqrt{\rm\sf M}x}\\
%& \Longleftrightarrow \norm{}{\sqrt{\rm\sf m}{C}^{-1/2}y} \le  \norm{}{{L}^{1/2}y} \le \norm{}{\sqrt{\rm\sf M}C^{-1/2}y} \\
%& \Longleftrightarrow {\rm\sf m}\, C^{-1} \le L \le  {\rm\sf M}\, C^{-1}.
%\end{align*}
%(2) The inequalities on the left imply that $L\in GL^{(+)}(\Hil)$. Then the result follows from Prop. \ref{sec:controlframclass2}.
%\eofproof 
\noi This result shows the equivalence of Eqs.(3.11) and (3.12) of Ref.\cite{bogdvan1} under the given conditions, which was stated there without proof.

\vspace*{-2mm}
\subsection{Numerical aspects of controlled frames}

As a short remark, we note that $(C \psi_n)$ need \emph{not} be related to a dual frame, except in the case $M = m$, when $C = L^{-1}$ by necessity. 
But, for finding the canonical dual frame algorithmically from \eqref{sec:contrframineq1}, 
we know that $L_C$ is invertible. In particular this means 
that $L_C^{-1} C = \left( C L \right)^{-1} C = L^{-1}$. So finding a $C$ such that $\Psi$ forms a controlled frame, with nice numerical properties, 
is equivalent to preconditioning.\cite{trebau1}  
This was the main motivation for introducing controlled frames in Ref.\cite{bogdvan1}.  
For this to be effective, the bounds of the operator, which are clearly also bounds for the spectrum of the operator, should be close to each other. 
%A clustered spectrum not only speeds up the convergence speed of iterative schemes, but also keeps the condition number low.
%So how are ${\rm\sf m}_L,{\rm\sf M}_L$ and ${\rm\sf m}_{CL}, {\rm\sf M}_{CL}$ related? 
It is straightforward to show:
\begin{cor.} Let $C$ be a self-adjoint operator and let $\Psi$ be a $C$-controlled frame. 
Denote by $({\rm\sf m}_{CL}, {\rm\sf M}_{CL})$,  $({\rm\sf m},{\rm\sf M})$ and  $({\rm\sf m}_C,{\rm\sf M}_C)$  any bounds for
the controlled frame operator $L_C$, the frame operator $L$, and the operator $C$,  respectively.
 Then, 

\qquad\begin{minipage}{12cm}\bei
\item[(i)] ${\rm\sf m'}  = \dis\frac{{\rm\sf m}_{CL}}{{\rm\sf M}_C}, \;{\rm\sf M'} = \dis\frac{{\rm\sf M}_{CL}}{{\rm\sf m}_C} \;$
 \mbox{are bounds for $L$}; 
\\[1mm]
\item[$\;$(ii)]  ${\rm\sf m}'_{C}  = \dis\frac{{\rm\sf m}_{CL}}{\rm\sf M}, \; {\rm\sf M}'_{C} = \dis\frac{{\rm\sf M}_{CL}}{\rm\sf m} \; $
 \mbox{are bounds for $C$}; 
\\[1mm]
\item[\hspace*{5mm}(iii)]  ${\rm\sf m}'_{CL}  = {\rm\sf m} {\rm\sf m}_C, \; {\rm\sf M}'_{CL} = {\rm\sf M} {\rm\sf M}_C \;$ 
 \mbox{are bounds for $L_{C}$}. 
\eni
\end{minipage}

\noi If two bounds are optimal in the above equations, the resulting third one is optimal, too.
\end{cor.}
%\noi{\bf Proof.}  Corollary \ref{sec:commutposoperbound1} gives us the equalities.
%Now suppose that two bounds are optimal (e.g. ${\rm\sf m}$ and ${\rm\sf m}_C$) and the third one (e.g. ${\rm\sf m}_{CL}$) is not, 
%so for example there is a $\widetilde{{\rm\sf m}}_{CL} > {{\rm\sf m}}_{CL}$. But this means that 
%$
%\widetilde{{\rm\sf m}} := \frac{\widetilde{{\rm\sf m}}_{CL}}{{\rm\sf M}_C} > \frac{{{\rm\sf m}}_{CL}}{{\rm\sf M}_C} = {\rm\sf m}
% $
%would also be a lower bound for $L$. This   contradicts the optimality of ${\rm\sf m}$.
%\eofproof
%\vspace{2mm}

This means that, if we find a $C$  such that ${\rm\sf m}_{CL} \cong {\rm\sf M}_{CL}$, we get a very efficient scheme, in the sense that:
\bei
\item[{\bf .}] The resulting algorithm is much more stable, according to the remarks made in Section \ref{sec:numiss0}. One has indeed
 $\kappa(L_C) \le\frac{{\rm\sf M}_{CL}}{{\rm\sf m}_{CL}}$.

\item[{\bf .}] Let $\epsilon :=\frac{{\rm\sf M}_{CL} - {\rm\sf m}_{CL} }{{\rm\sf M}_{CL} + {\rm\sf m}_{CL}}$. 
Using a Neumann algorithm, we get a good approximation of the inverse operator 
already in the first iteration $g_1$. Indeed, $\norm{}{f -g_1} \le \epsilon \norm{\CC^n}{f}$.
\eni
So, as stated in Ref.\cite{bogdvan1}, although controlled frames and ``standard" frames are mathematically equivalent, these different `viewpoints' of frames give 
opportunities for efficient implementations. For general frames, it seems difficult to find an appropriate preconditioning matrix, but for wavelet frames this technique
 is used in the above-mentioned paper.\cite{bogdvan1} For Gabor frames, a way to find advantageous preconditioning matrices is presented in Ref.\cite{xxlfei1}.

\section{Weighted frames} \label{sec:weightfram0}
\vspace*{-2mm}
\bedefi 
Let $\Psi = (\psi_n : n \in \Gamma)$ be a sequence of elements in $\Hil$ and   $( w_n : n \in \Gamma ) \subseteq \RR^+$ a sequence of  positive weights. 
This pair is called   
a \emph{w-frame} of $\Hil$ if there exist constants ${\rm\sf m}>0 $ and ${\rm\sf M} < \infty$ such that
\begin{equation} \label{eq:frame}
{\rm\sf m} \norm{}{f}^2 \le \sum \limits_{n\in \Gamma}  w_n \left| \left< f , \psi_n \right> \right|^2 \le {\rm\sf M} \norm{}{f}^2 .
\end{equation} 
Alternatively, given a sequence of complex numbers  $(\omega_n) \subseteq \CC$, we call $\Psi = (\psi_n) $  a 
 \emph{weighted frame}  if the sequence  $(\omega_n\, \psi_n)$ is a frame.\footnote{The two terms `weighted frame'
and `w-frame' were used interchangeably in Ref.\cite{bogdvan1}; here we make a difference.}
 
\end{def.}
The two definitions are clearly equivalent, by putting  $w_n = \left| \omega_n \right|^2$, resp., $\omega_n = \sqrt{w_n} \epsilon_n$, 
where $\epsilon_n \in \CC$ with $\left| \epsilon_n \right| = 1$. 

Weighted frames are related to {\em signed frames}.\cite{pewa02}  The latter are Bessel sequences coupled with weights $(w_n = \pm 1)$, 
that fulfill an inequality similar to \eqref{eq:frame}. 
Signed frames are  equivalent to w-frames for which negative weights are allowed and the weights are bounded.  
In this paper, the sign of the weight is either not included or not significant for the definitions given above, 
and the main results in  Ref.\cite{pewa02} also differ in focus from the ones given here.

\vspace*{-2mm}
\subsection{w-frames as controlled frames} \label{sec:weightcontr0}

It is clear that, if an operator is diagonal on a given sequence and together they form a controlled frame, then the concept of weighted frame is just a 
special case of controlled frames. 
But we can show that we cannot get all possible cases of  
\emph{w-frames}   in that way and so we cannot apply the result in Section \ref{sec:contrfram1} 
to the general w-frame case. 

\begin{prop.}\label{sec:contrweightsemi1} Let $C \in GL(\Hil)$ be self-adjoint and diagonal on $\Psi= (\psi_n)$ and assume it generates a controlled frame. 
Then the sequence $(w_n)$, which verifies the relations  $C \psi_n = w_n \psi_n$, is semi-normalized \footnote{A sequence $(c_n)$ is called semi-normalized  
if there are bounds $b \ge a > 0$, such that $a \le | c_n | \le b$.} and positive. Furthermore $C = \MM_{{\bf w},\tilde \Psi, \Psi}$.
%
%If $w_n \not= w_k$ for all $k,n$ then $C = I$.
\end{prop.}
\begin{proof}
%
%Let us suppose the operator $C \in GL(\Hil)$ is self-adjoint, diagonal on $\Psi$, i.e.,  $C \psi_n = w_n \psi_n$, and turns $\Psi$  into a controlled frame. 
By Proposition \ref{sec:posinvertop1}, we get the following result for $C^{1/2}$:
$$
{\rm\sf m}_C \norm{}{f}^2 \le \norm{}{{C}^{1/2}f}^2 \le {\rm\sf M}_C \norm{}{f}^2.
$$
As $C \psi_n = w_n \psi_n$, clearly ${C}^{1/2} \psi_n = \sqrt{w_n} \psi_n$. 
Applying the inequalities above to the elements of the sequence, we get
$
 0 < {\rm\sf m}_C \le w_n \le {\rm\sf M}_C
$.
%
% Let us suppose as in Prop. \ref{sec:contrweightsemi1} that $C$ is self-adjoint, bounded, forms a controlled frame with $\Psi$ and acts diagonally on $\Psi$ with the weights ${\bf w}=(w_k)$.
 
% It is easy to show, using for example the matrix representation of operators using frames \cite{xxlframoper1}, that 
 Clearly, the only possible operator $C$ that could fulfill the conditions would be the multiplier 
 $$C f = \sum \limits_k \left< f , \tilde \psi_k \right> w_k \psi_k = \MM_{{\bf w},\tilde \Psi, \Psi} f .$$ 
%
%For weights $w_k \not= w_l$, $\psi_k$ and $\psi_l$ are orthogonal to each other. According to Prop. \ref{sec:controlframclass2},
% $L$ and $C$ commute. In particular they have the same eigenspaces and so $L$ and also $L^{-1}$ act diagonally on $\Psi$. So let $L^{-1} \psi_k = w_k \psi_k$. Then the dual frame of $\Psi$ is just the weighted frame $(w_k \psi_k)$. 
%  Thus $\MM_{{\bf w},\tilde \Psi, \Psi}  = I $.
\end{proof}

The  relation  between controlled frames and weighted frames   for non-self-adjoint operators $C$ is not obvious.
The reason is that, for nonexact frames, a definition by $U \psi_n = w_n \psi_n$ 
is not applicable.\cite{xxlframoper1} 
%Indeed it is not well-defined, since the vectors  $\psi_n$ are not linearly independent in general. If this shortcoming is circumvented by setting $U f = \sum \limits_n \left< f , \widetilde \psi_n  \right> w_n \psi_n$,
% then in general $U \psi_n \neq w_n \psi_n$. See \cite{xxlframoper1} for more on this.
% \\
%
%For the special case of regular Gabor frames (see also Section \ref{sec:gabornum0}), this is only possible for Gabor frames that are already tight and 
%where the weights are constant.\cite{doerf1}  So, in particular:
%\begin{cor.} Let $G = \left( \pi ( \lambda) g \right)$ be a regular Gabor frame and $(w_\lambda)$ weights that are not all equal. 
%Then the pair $\left( \left(w_\lambda\right),G \right)$ cannot be described 
%as a controlled frame.
%\end{cor.}
\vspace*{-2mm}
\subsection{Semi-normalized weights}

As a converse to the first part of Prop. \ref{sec:contrweightsemi1},
  a frame weighted by a semi-normalized sequence  is always a frame. Indeed, 
\begin{lem.} \label{sec:framepertmul1} Let $( \omega_n )$ be a semi-normalized sequence with bounds $a$,$b$. 
Then if $\left( \psi_n \right)$ is a 
frame with bounds ${\rm\sf m}$ and ${\rm\sf M}$,  then $\left( \omega_n \psi_n \right)$ is also a frame with bounds $a^2 {\rm\sf m}$ and $b^2 {\rm\sf M}$. 
 The sequence $\left({\omega_n}^{-1} \tilde \psi_n \right)$ is a dual frame of $( \omega_n \psi_n )$.
\end{lem.}
\noi{\bf Proof.} 
Since $\left| \left< f, \omega_n \psi_n \right> \right|^2 = \left| \omega_n \right|^2 \left| \left< f, \psi_n \right> \right|^2 $, we get
$$
\Delta:=\sum \limits_n \left| \left< f, \omega_n \psi_n \right> \right|^2 = \sum \limits_n \left| \omega_n \right|^2 \left| \left< f, \psi_n \right> \right|^2.   
$$
Thus
$
\Delta\le b^2 \sum \limits_n \left| \left< f, \psi_n \right> \right|^2 \le b^2 {\rm\sf M} \norm{}{f}^2
 $. 
In addition, 
$$
 \Delta \ge a^2 \sum \limits_n \left| \left< f, \psi_n \right> \right|^2 \ge a^2 {\rm\sf m} \norm{}{f}^2 .$$ 
As $\sum \limits_n \left< f , \omega_n \psi_n \right> {\omega_n}^{-1} \tilde \psi_n = \sum \limits_n \left< f , \psi_n \right> \tilde \psi_n = f$, 
these two sequences are dual. Since ${\omega_n}^{-1} $ is bounded, $\left( {\omega_n}^{-1}  \tilde \psi_n \right)$ is a Bessel sequence dual to a frame. 
Therefore,\cite{ole1}  it is a dual frame of $(\omega_n \psi_n)$. 
\vspace*{-1cm} 
\eofproof
\vspace*{1cm} 

\noi{\bf Remarks}
\begin{enumerate}
\item The weighted dual frame is a dual, but not the canonical dual. As an example, consider the Parseval frame, i.e.,  self-dual frame, 
$$\Psi = \left\{ \left( \begin{array}{c} 
\frac{2}{\sqrt{6}}  \\
0      
\end{array}
\right), \left( \begin{array}{c} 
-\frac{1}{\sqrt{6}}  \\
 \frac{1}{\sqrt{2}} 
\end{array}
\right), \left( \begin{array}{c} 
 - \frac{1}{\sqrt{6}} \\
 -\frac{1}{\sqrt{2}}
\end{array}
\right) 
\right\}$$
 with weights $(\omega_1, \omega_2, \omega_3) = (\frac{1}{2}, 1 , 2)$. Following Lemma \ref{sec:framepertmul1} 
  $({\omega_k}^{-1}  \psi_k)$ forms a dual frame,
 but it is not identical to the canonical dual frame (see Figure \ref{fig:weightdual1}).
The relationship between these two duals will be investigated in Section \ref{sec:gabornum0} for the case of Gabor frames.
\begin{figure}[t]
\begin{center}
\includegraphics[width=0.4\textwidth]{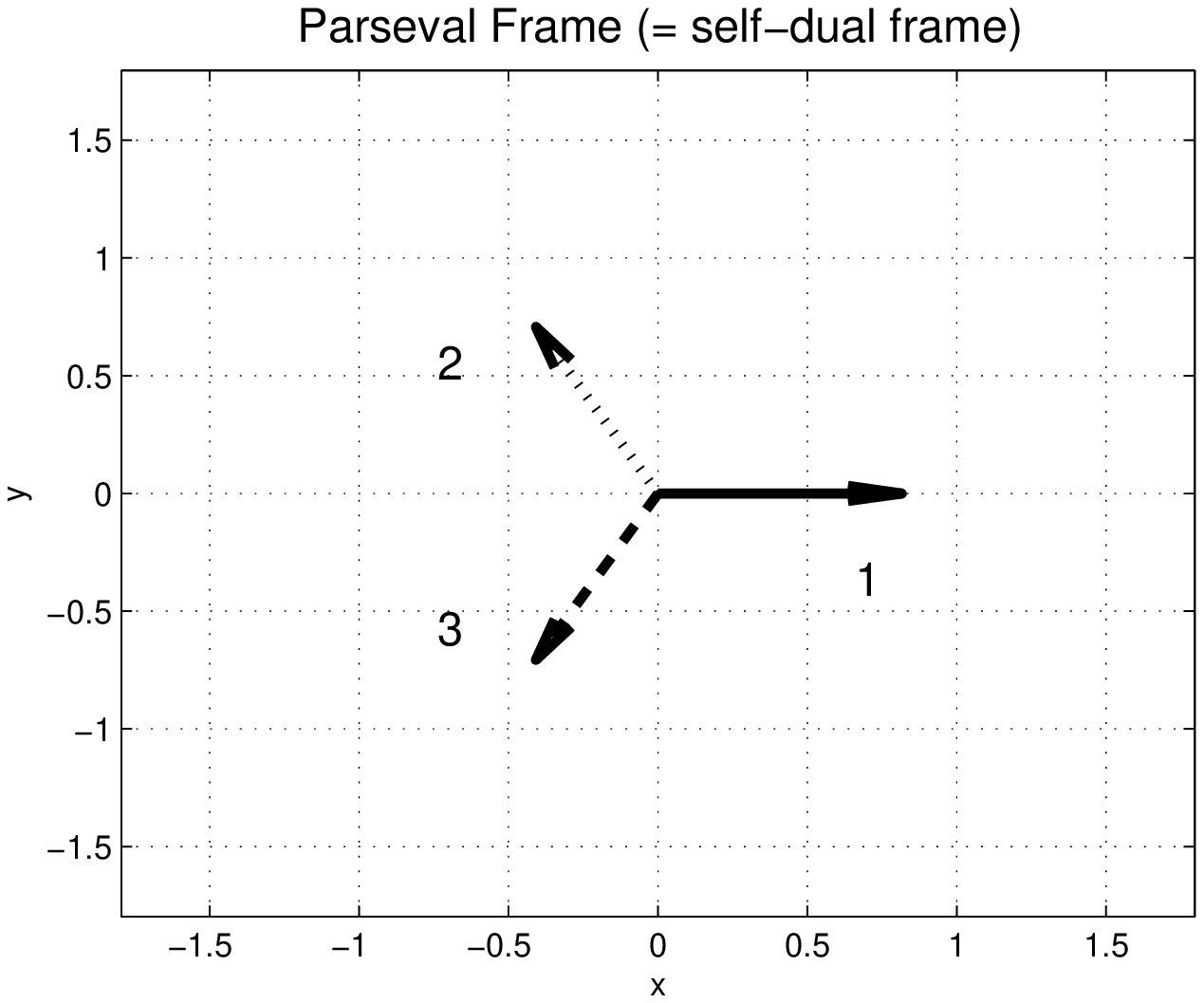} \hspace{3mm}
\includegraphics[width=0.4\textwidth]{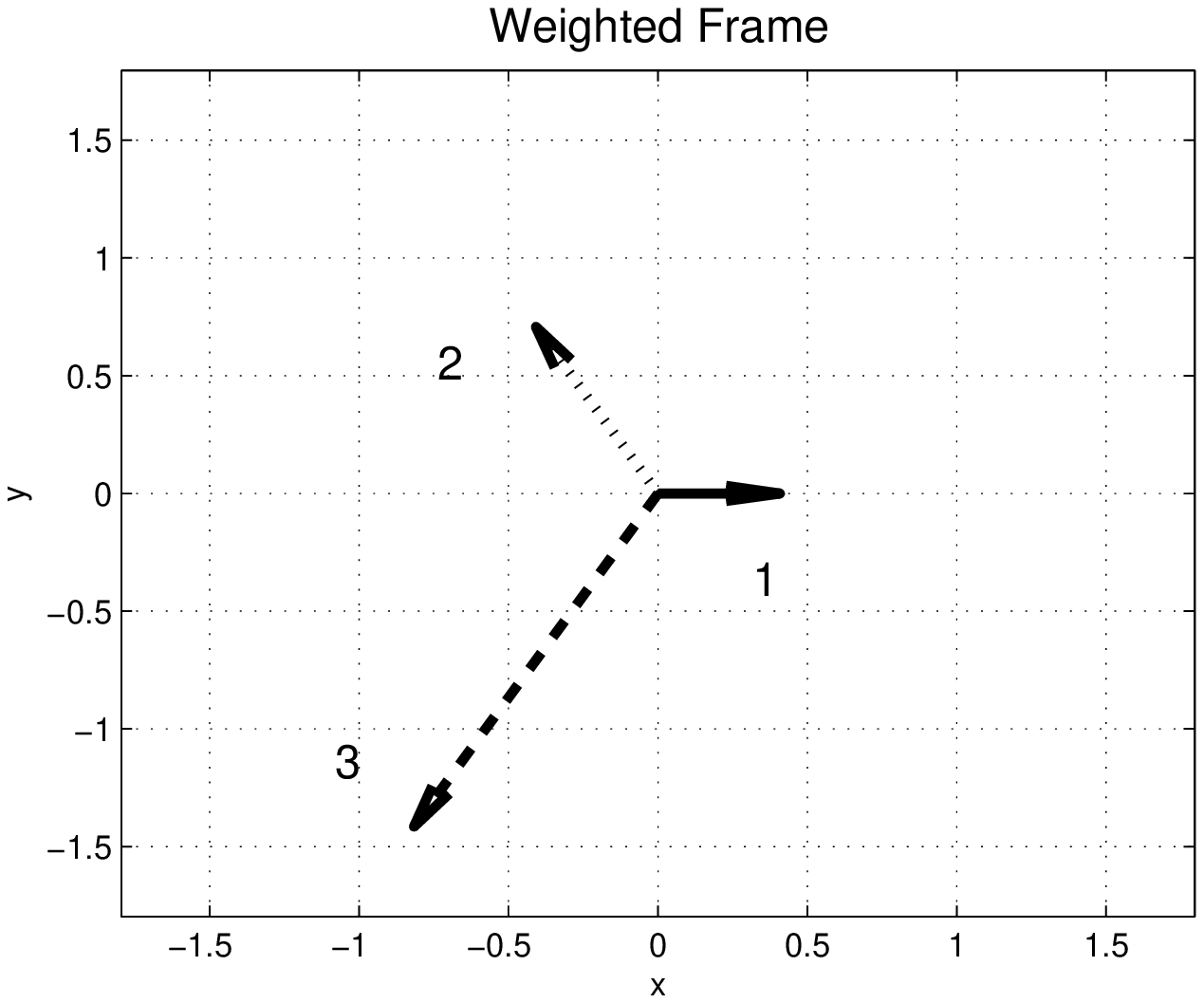}

\includegraphics[width=0.4\textwidth]{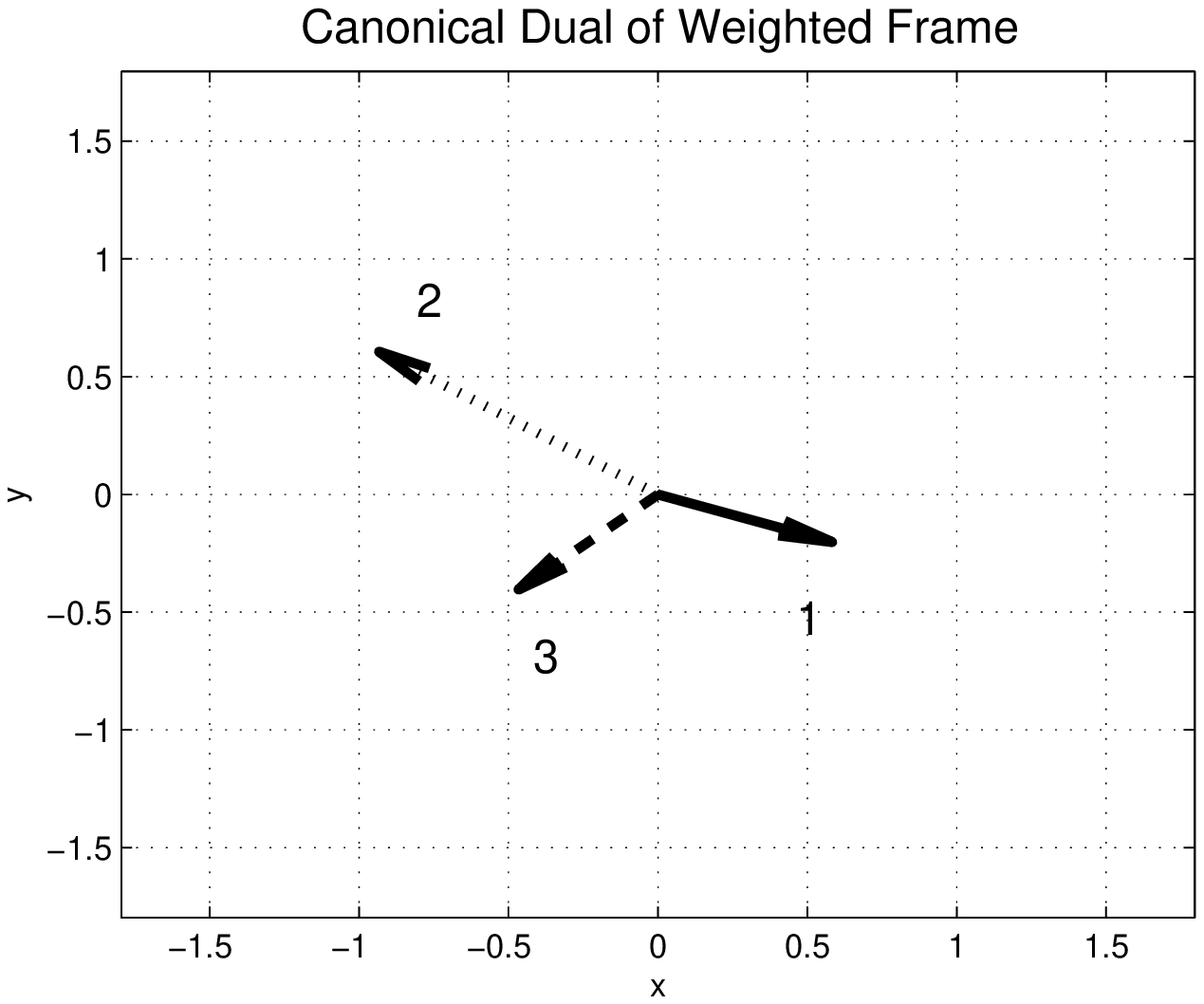} \hspace{3mm}
\includegraphics[width=0.4\textwidth]{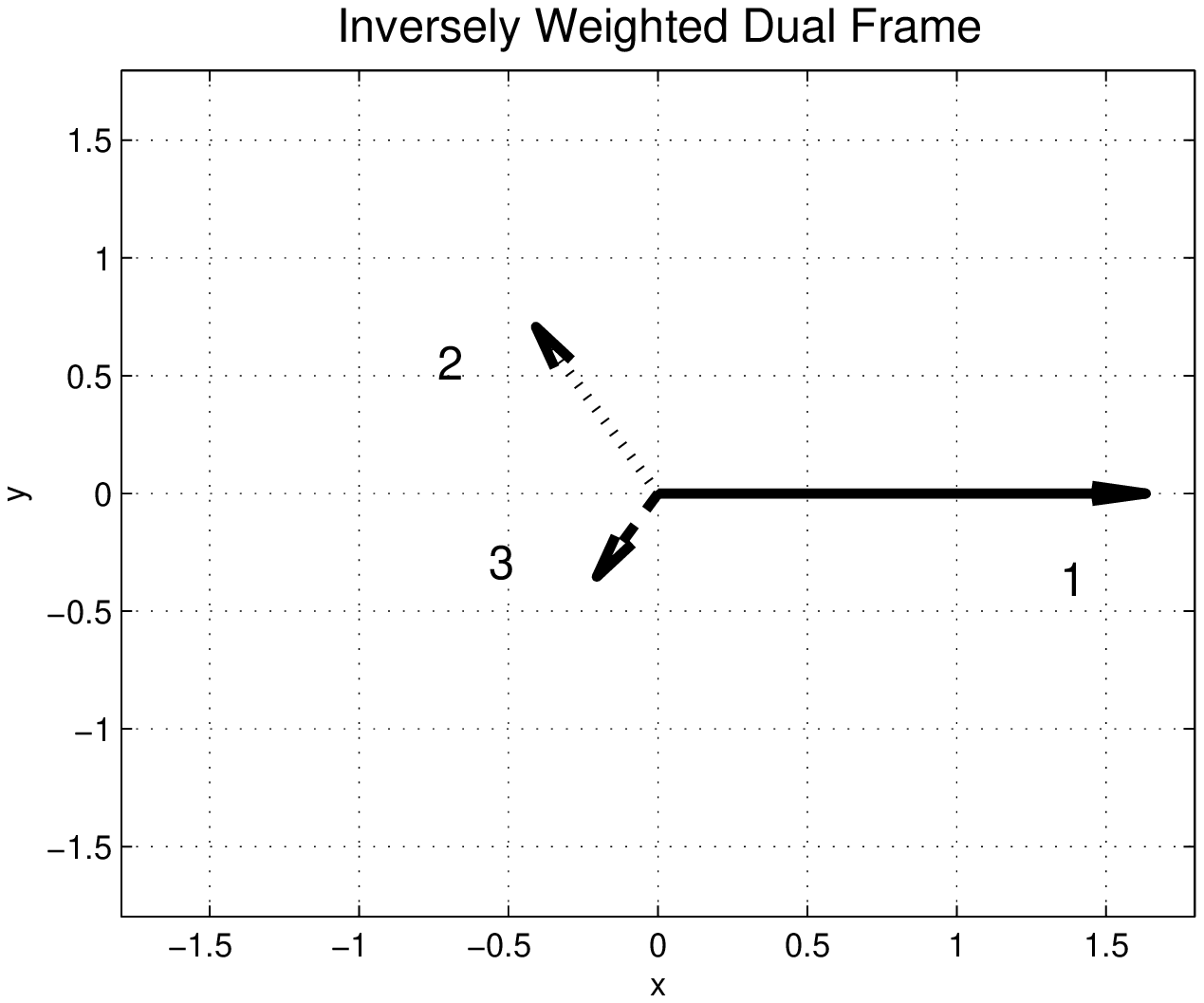}

\end{center}
\caption{\label{fig:weightdual1}\small \em Comparing the canonical dual of a weighted frame to the reciprocal weighted dual frame. 
(Top left:) original frame. (Top right:) Weighted frame. (Bottom left:) Canonical dual of weighted frame. 
(Bottom right:) The canonical dual of the original frame, i.e.,  the frame itself, weighted by the inverse weights.} 
\end{figure} 

\item Frames with weights $(w_n = \pm 1)$ fulfill the conditions of this Lemma, so they are always a frame with the same bounds. 
A dual frame is obtained just by applying the weights on the dual of the frame. This means that for a signed frame\cite{pewa02} that is also a frame, 
the dual can be easily calculated.  

\item The condition for semi-normalized sequences in Lemma \ref{sec:framepertmul1} is necessary. 
It is in general not enough for the weights to be strictly positive,   $w_n > 0$, for all $n$.  To give an example, let $(e_n)$ be an orthonormal basis 
in $\Hil$ with index set $\mathbb{N}$ and $\psi_n = \frac{1}{n} e_n$. This is not a frame, since this sequence does not fulfill the lower frame condition. 

\item Using the sequence $(\psi_n)$ above  with the weights $w_n = n$ also shows that in general a weighted frame need not be a frame, 
if the semi-normalized condition is not fulfilled. Furthermore this shows that Lemma \ref{sec:framepertmul1} is not reversible.
 There are cases where weights that are not semi-normalized lead to weighted frames. There are even cases where unbounded sequences lead to weighted frames. 
\end{enumerate}
\vspace*{-2mm}
\subsection{Connection to frame multipliers} \label{sec:frammul0}
 
The concept of  weighted frames  is connected to that of  frame multipliers.\cite{xxlmult1}  
%Let us investigate this connection in more detail.
\begin{lem.} \label{sec:framepertmul2} 
Let $\Psi=\left( \psi_n \right)$ be  a frame for $\Hil$. Let ${\bf m}=( m_n )$ be a positive, 
semi-normalized sequence.% with bounds $a,b$. % such that $0 < a \le \left| m_n  \right| \le b \ \forall k$. 
Then the multiplier $\MM_{{\bf m}, \Psi}$ is  the frame operator of the frame $( \sqrt{m_n} \psi_n )$ and therefore it  is positive, self-adjoint and invertible. 
\\
If $(m_n)$ is negative and semi-normalized, then  $\MM_{{\bf m}, \Psi}$ is negative, self-adjoint and invertible.
\end{lem.}
\noi{\bf Proof.} 
$$ 
\MM_{{\bf m}, \Psi}f  = \sum \limits_n m_n \left< f, \psi_n \right> \psi_n = \sum \limits_n \left< f, \sqrt{m_n} \psi_n \right> \sqrt{m_n} \psi_n.
$$
By Lemma \ref{sec:framepertmul1}, $( \sqrt{m_n} \psi_n )$ is a frame. Therefore $M_{{\bf m}, \Psi}= L_{( \sqrt{m_n} \psi_n )}$ is positive and invertible.

Let $m_n <0$ for all $n$, then $m_n = - \sqrt{\left| m_n \right|}^2$. Therefore 
$$
 \MM_{{\bf m}, \Psi} = - \sum \limits_n \left< f, \sqrt{\left| m_n \right|} \psi_n \right> \sqrt{\left| m_n \right|} \psi_n 
= - L_{( \sqrt{\left| m_n \right|} \psi_n)} .
\vspace{-8mm}$$
\eofproof

%\vspace{8mm}
This can be extended to 
\begin{theorem} Let $(\psi_n)$ be a sequence of elements in $\Hil$. Let $(w_n)$ be a sequence of positive, semi-normalized weights. 
Then the following properties are equivalent:
\begin{enumerate}
\item $(\psi_n)$ is a frame.
\item $\MM_{{\bf m}, \Psi}$ is a positive and invertible operator.
\item There are constants ${\rm\sf m}>0 $ and ${\rm\sf M} < \infty$  such that
$$
 {\rm\sf m} \norm{}{f}^2 \le \sum \limits_{n \in \Gamma} w_n \left| \left< f , \psi_n \right> \right|^2  \le {\rm\sf M} \norm{}{f}^2,
 $$
i.e., the pair $(w_n)$,$(\psi_n)$ forms a w-frame.
\item $(\sqrt{w_n} \psi_n)$ is a frame.
\item $\MM_{{\bf w'}, \Psi'}$ is a positive and invertible operator for any positive, semi-normalized sequence $(w'_n)$.
\item $(w_n \psi_n)$ is a frame, i.e., the pair $(w_n)$,$(\psi_n)$ forms a weighted frame.
\end{enumerate}
\end{theorem}
\noi{\bf Proof.}  By Lemma \ref{sec:framepertmul2}, (1) $ \Longrightarrow $ (2). By Proposition \ref{sec:posinvertop1}, (2) $\Longleftrightarrow$ (3). 
By the definition of w- and weighted frames we get $(3) \Longleftrightarrow (4)$.

If $(w_n)$ is positive and semi-normalized, so is the sequence $({w_n}^{-1})$. With Lemma \ref{sec:framepertmul2} and application of (1) $
\Longrightarrow$ (4),  we get  (4) $\Longrightarrow$ (1).

Use the sequence $(w'_n)$ in the above argument to get (1) $\Longleftrightarrow$ (5).

Finally,   
$(w_n^2)$ is also a positive, semi-normalized sequence, therefore the results above  show the equivalence of (6) with the rest.
\eofproof

Clearly this could also be easily extended to negative weights.
\vspace*{-2mm}
\section{Numerical results for general frames} 
\label{sec:weightnum0}

Now the practical question is obviously, how can one find weights   such that a weighted frame becomes `as tight as possible', such that the quotient of the bounds, i.e.,  the condition number of the frame, becomes smaller? 
This would give a way to calculate a dual in a more efficient way, although one does not obtain the canonical dual in general.
In the sequel, we consider several possibilities, namely, $\ell^p$-type weights and weights obtained through approximation by a frame multiplier,
and we study the performance of each type of weight (in finite dimensions, of course).
 Many numerical tests have been performed, but, for the sake of conciseness, we present only the most significant results.
\vspace*{-2mm}
\subsection{First method: $\ell^p$-weights }

We are looking for a measure of how important a single frame element is, how dependent it is on the other elements. 
With this in mind,  given a frame $\Psi = \{\psi_k,\, k=1,\ldots,M\} \;(M\le \infty)$, put 
$$
\omega^{(2)}_n = \frac{\left\| \psi_n \right\|}{\sqrt{\sum \limits_k \left| \left< \psi_n , \psi_k \right> \right|^2}}.
$$ 
This weight can be motivated as a control of the importance of the side-diagonals of the Gram matrix, ${G_{\Psi}}_{k,l} = \left< \psi_k, \psi_l\right>$, 
by comparing a diagonal entry to the sum of the squares of the other entries on the same line. This is reminiscent also of the generalized 
Welch bound found in Ref.\cite{Waldron03}.
 For an orthogonal basis, the best weights should be given by the normalization, which is indeed achieved by this weight, since in this case 
one has $\omega^{(2)}_n = \frac{\left\| \psi_n \right\|}{\left\| \psi_n \right\|^2}$.

In order to measure the influence of the power chosen  in the definition of  $\omega^{(2)}_n$, we have also investigated
$$
\omega^{(p)}_n = \frac{\left\| \psi_n \right\|}{\left(\sum_k \left| \left< \psi_n , \psi_k \right> \right|^p\right)^{1/p}}
$$
for $p=4$ and $p=6$. In other numerical experiences, it could be observed that $p=1,3,5$ are not good choices compared to those ones. 

Finally we also  test the following weight:
$$
\omega^{(\infty)}_n = \frac{\norm{}{\psi_n}} {\sup_k \left| \left< \psi_n , \psi_k \right> \right| } .
$$
\vspace*{-2mm}
\subsection{Second method: Weights by best approximation with a frame multiplier}

These questions can also be translated into the frame multiplier context: Can the identity be written as a frame multiplier? Can it be approximated?

It is possible to find the best approximation of operators (in finite-dimensional discrete spaces) using the Hilbert-Schmidt norm (see Ref.\cite{xxlframehs07} for an algorithm).
 The symbol of the best approximation is the weight $\omega^{\rm (mult)}$, defined as follows:
%$$
% \omega^{\rm (mult)} =\sqrt{  \left(G_{\Psi}^{(2)} \right)^\dagger \norm{\Hil}{\psi_k}^2  }, 
%$$
$$ \omega^{\rm (mult)}_{n} =\sqrt{ \dis\sum_{k=1}^M \left[ \left( G_{\Psi}^{(2)} \right)^\dagger\right]_{nk }
  \norm{\Hil}{\psi_k}^2  }, 
  $$
where 
$ G_{\Psi}^{(2)}$ is the matrix
 $ \big(G_{\Psi}^{(2)} \big)_{pq} = |\left< \psi_q , \psi_p \right>|^2$
% ${G_{\Psi}}_{p,q} = \left< \psi_q , \psi_p \right>$, 
%the norm and square in $\left| G_{\Ppsi} \right|^2$ are taken elementwise 
and  $^\dagger$ denotes the pseudo-inverse. % of the resulting matrix. 
The resulting  matrix acts on the sequence $(x_k)_k = \left( \left( \norm{\Hil}{\psi_k} \right)^2\right)_k$. 
This method corresponds to finding the weight such that the frame operator of the weighted frame is as similar  as possible  to the identity (in the Hilbert-Schmidt topology).
%\begin{figure}[t]
%\begin{center}
%\includegraphics[width=0.48\textwidth]{histo_43_10000_new.eps}
%\includegraphics[width=0.48\textwidth]{histo_63_10000_new.eps}
%\includegraphics[width=0.48\textwidth]{histo_93_10000_new.eps}
%\includegraphics[width=0.48\textwidth]{histo_123_10000_new.eps}
%\end{center}
%\caption{ \label{fig:testtight1}\small \em  Frames  in $d=3$ dimensions: Improvement of condition number by weights $\omega^{(2)}$ (= `2-norm'), $\omega^{(4)}$ (= `4-norm'),
% $\omega^{(6)}$ (= `6-norm'), $\omega^{(\infty)}$ (= `inf.-norm') and $\omega^{\rm (mult)}$ (= `Multiplier'). 
%Top left: Frame with $M = 4$ elements; top right: $M = 6$; bottom left: $M = 9$; bottom right: $M = 12$}
%\end{figure}
\vspace*{-2mm}
\subsection{Procedure and results}

In order to compare the efficiency of the various types of weights, we restrict ourselves to finite dimension $d<\infty$ and create
random frames   by finding $M$ random elements $(M > d)$ and checking whether they span the whole space. For algorithms see Ref.\cite{xxlfinfram1}.  
In each case,  we calculate the condition number of the frame operator of the weighted frame and compare 
the different weights. For each of them, we test     whether it improves the condition number and whether it is the best among the given options.
This is repeated 10.000 times. 

In the following graphs, the color grey corresponds to the cases where the weighted frame improves the condition number of the frame matrix and 
black if the option is the best of the  three  given weights. So the sum of the black bars gives the total percentage of weights that  improve the condition number.

%First we present results for dimensionality $d = 3$ and increasing redundancy, the number of elements are $M = 4,6,9$ and 12, respectively. 
%The result is shown in Figure \ref{fig:testtight1}. For this test, at least one of the given weights is improving the condition number in 95.17 
%\%, 99.14 \%, 100 \%, resp.100 \% of the cases.
%
% In Figure \ref{fig:testtight2} we present results for dimensionality $d = 10$ and increasing redundancy $(=M/d)$, the number of elements are chosen as $M = 11,20,30$ and $40$. 
%The condition number is improved in  $92.51 \%$, $99.25 \%$, $99.97 \%$, resp. $100 \%$ of the tests.

%\begin{figure}[t]
%\begin{center}
%\includegraphics[width=0.48\textwidth]{histo_1110_10000_new.eps}
%\includegraphics[width=0.48\textwidth]{histo_2010_10000_new.eps}
%
%\includegraphics[width=0.48\textwidth]{histo_3010_10000_new.eps}
%\includegraphics[width=0.48\textwidth]{histo_4010_10000_new.eps}
%\end{center}
%\caption{\label{fig:testtight2}\small \em  Frames  in $d=10$ dimensions, with the same conventions as in Figure \ref{fig:testtight1}.
% Top left: Frame with $M = 11$ elements; top right: $M = 20$; bottom left: $M = 30$; bottom right: $M = 40$.} 
%\end{figure}

\begin{figure}[t]
\begin{center}
\includegraphics[width=0.4\textwidth]{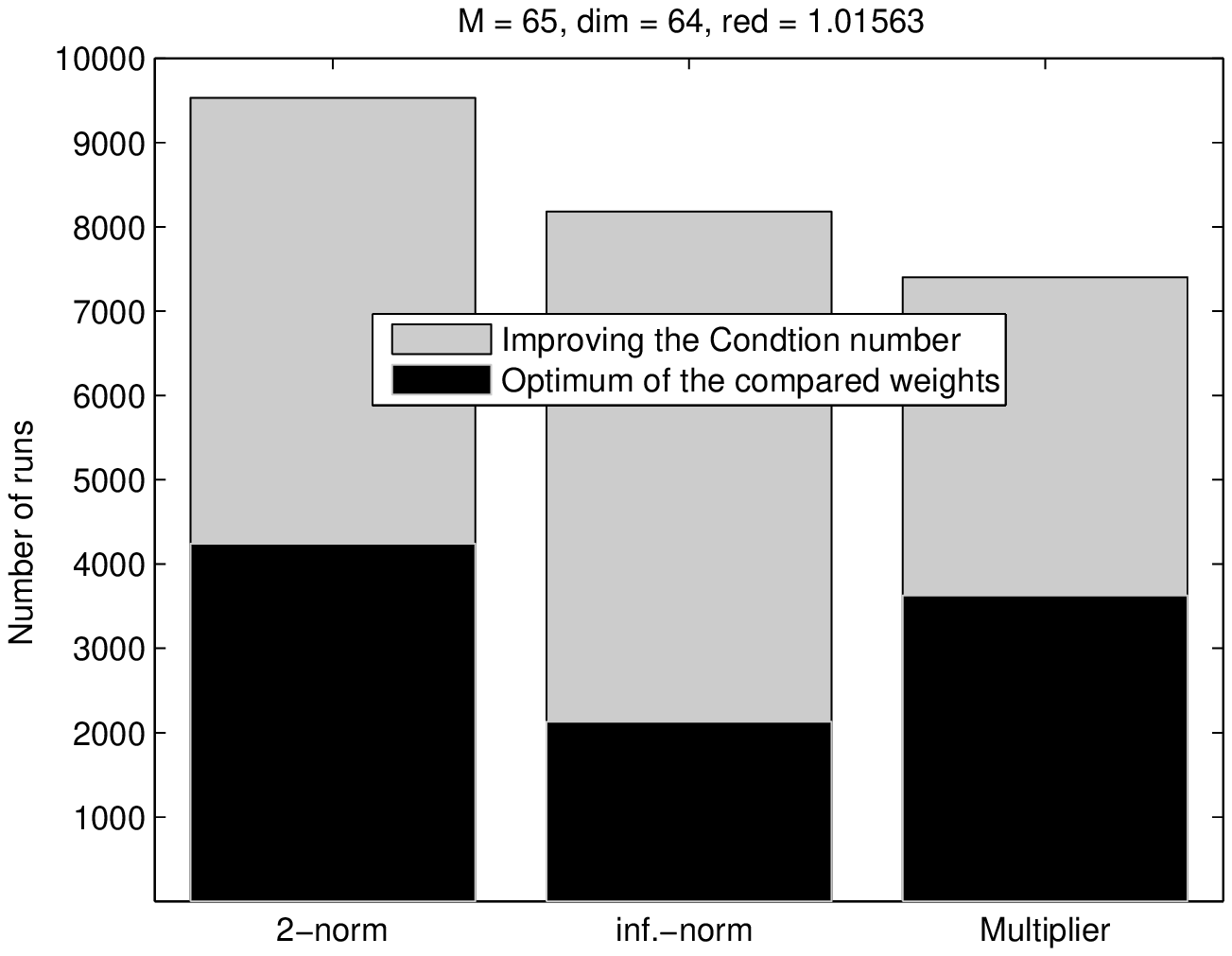} \vspace{3mm}
\includegraphics[width=0.4\textwidth]{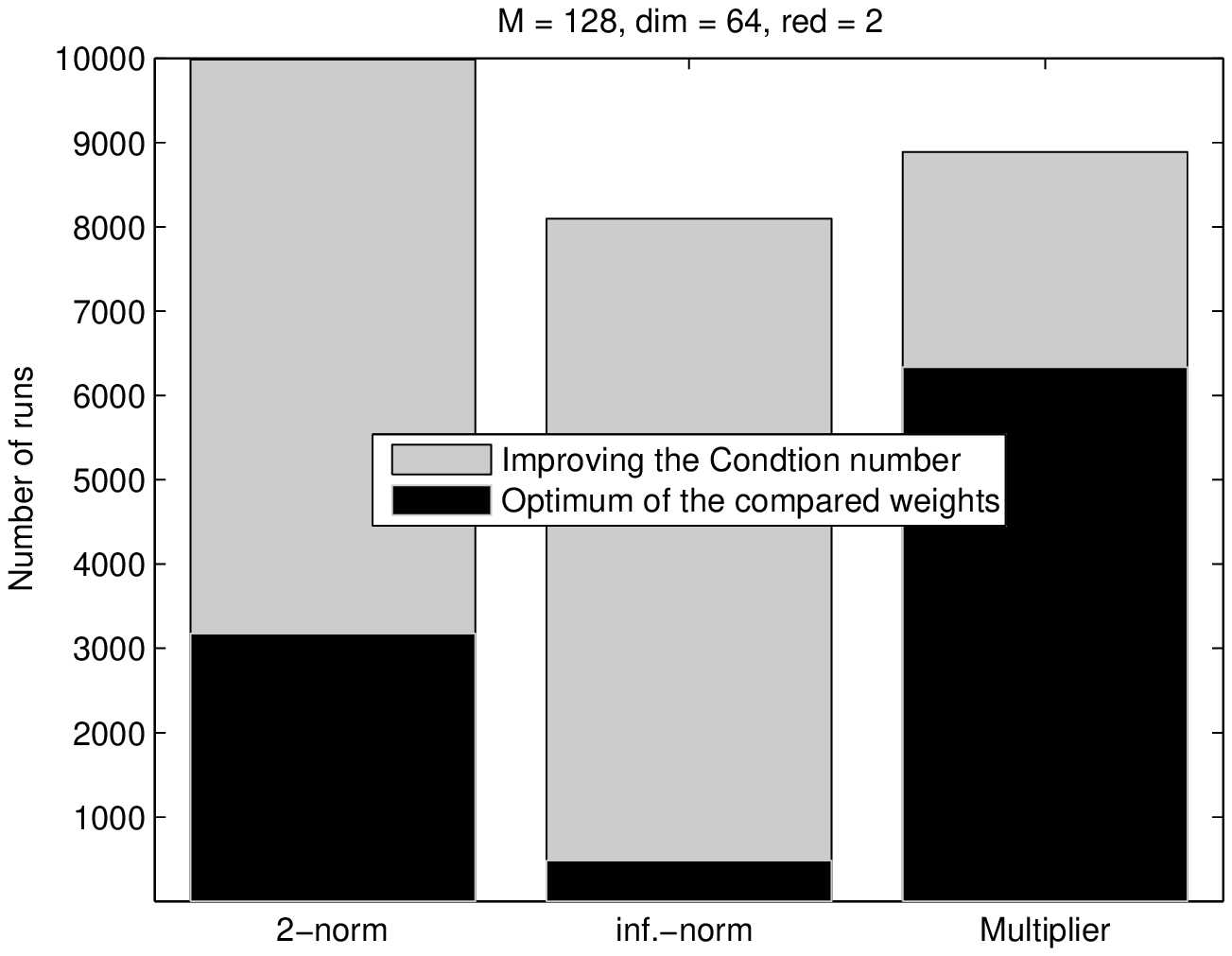}

\includegraphics[width=0.4\textwidth]{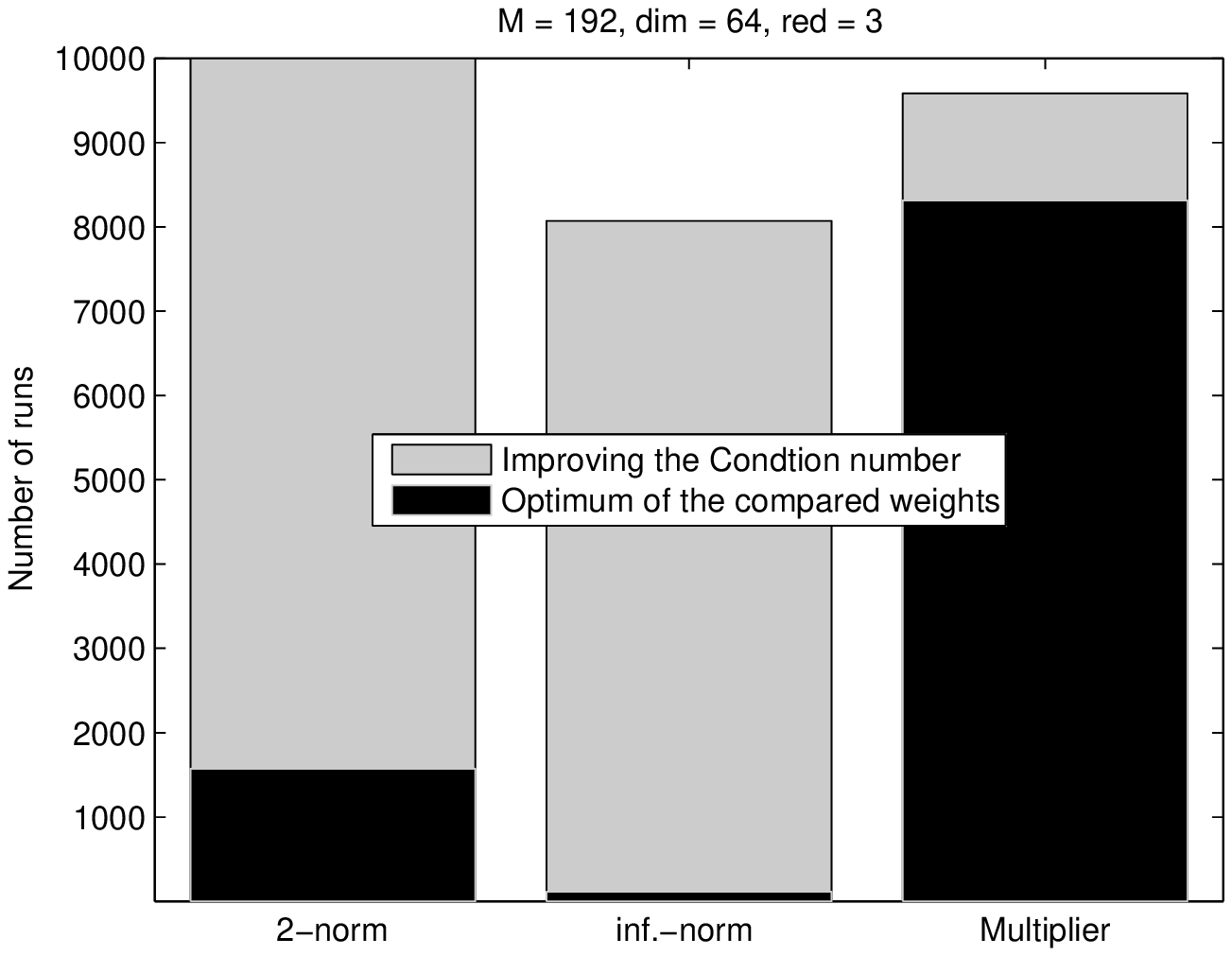}\vspace{3mm}
\includegraphics[width=0.4\textwidth]{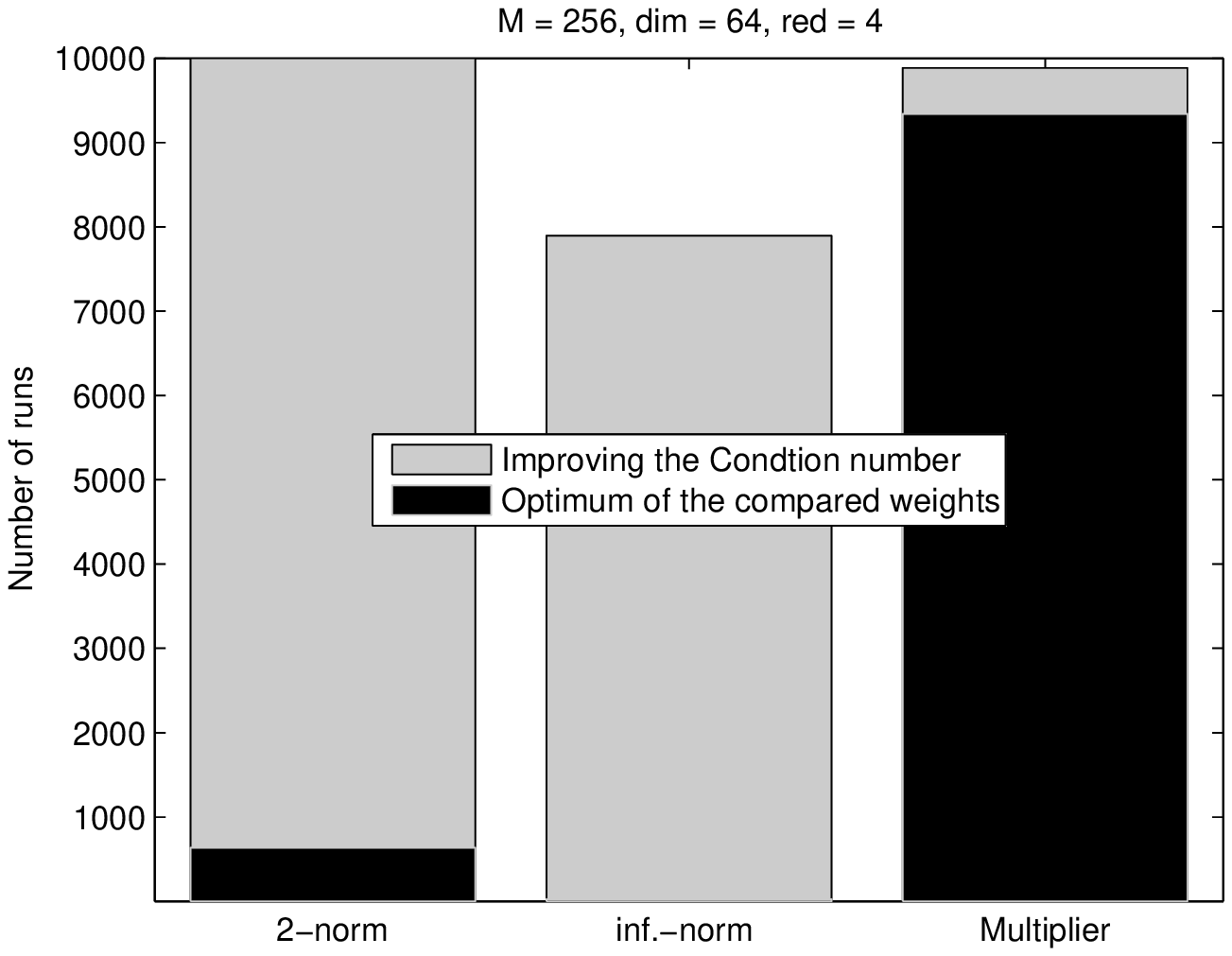}
\end{center}
\caption{\label{fig:testtight3}\small \em  Frames  in $d=64$ dimensions. Improvement of condition number by weights $\omega^{(2)}$ (= `2-norm'),
, $\omega^{(\infty)}$ (= `inf.-norm') and $\omega^{\rm (mult)}$ (= `Multiplier').
Top left: Frame with $M = 65$ elements; top right: $M = 128$; bottom left: $M = 192$, bottom right: $M = 256$.} 
\end{figure}

\begin{figure}[t]
\begin{center}
\includegraphics[width=0.4\textwidth]{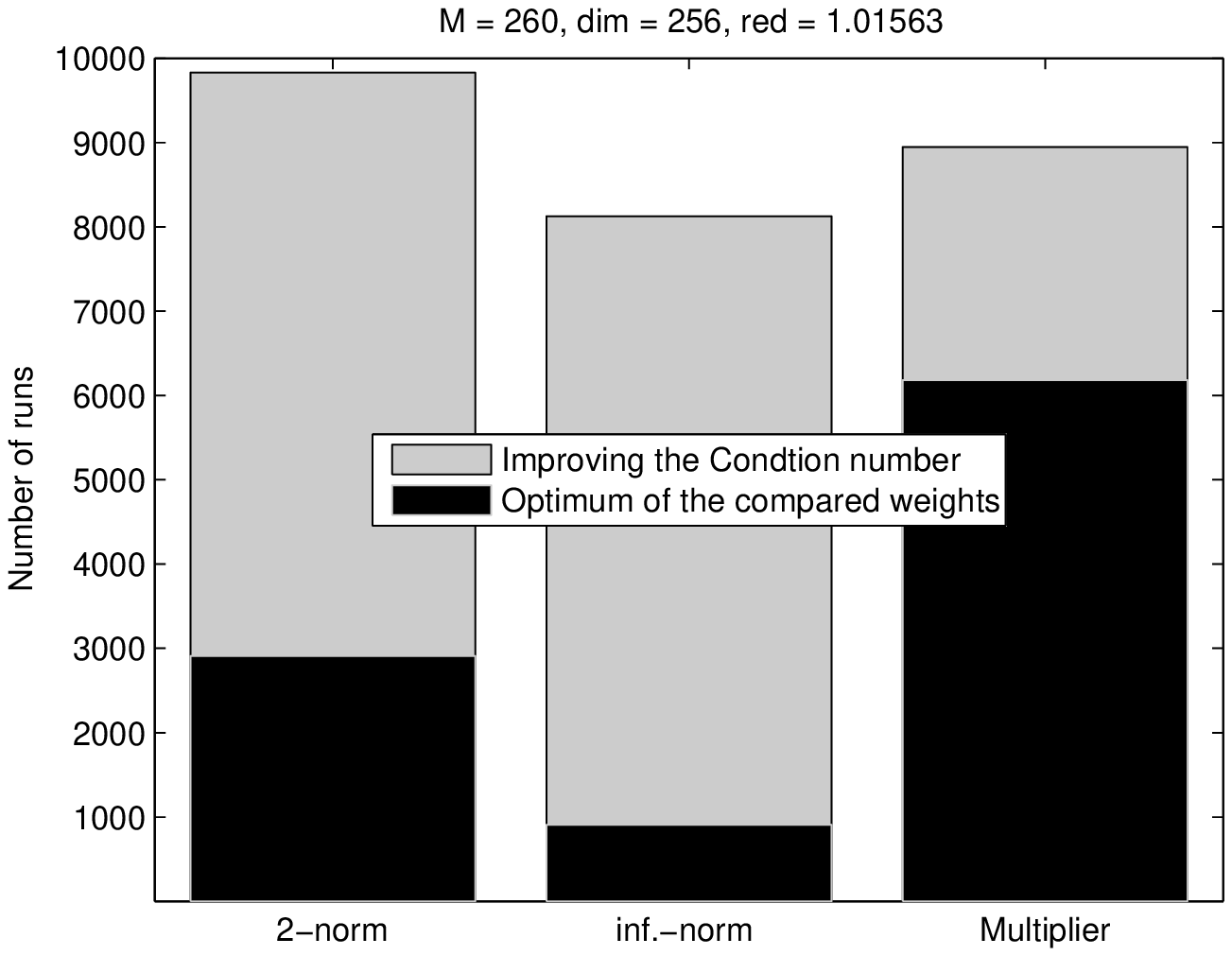}\vspace{3mm}
\includegraphics[width=0.4\textwidth]{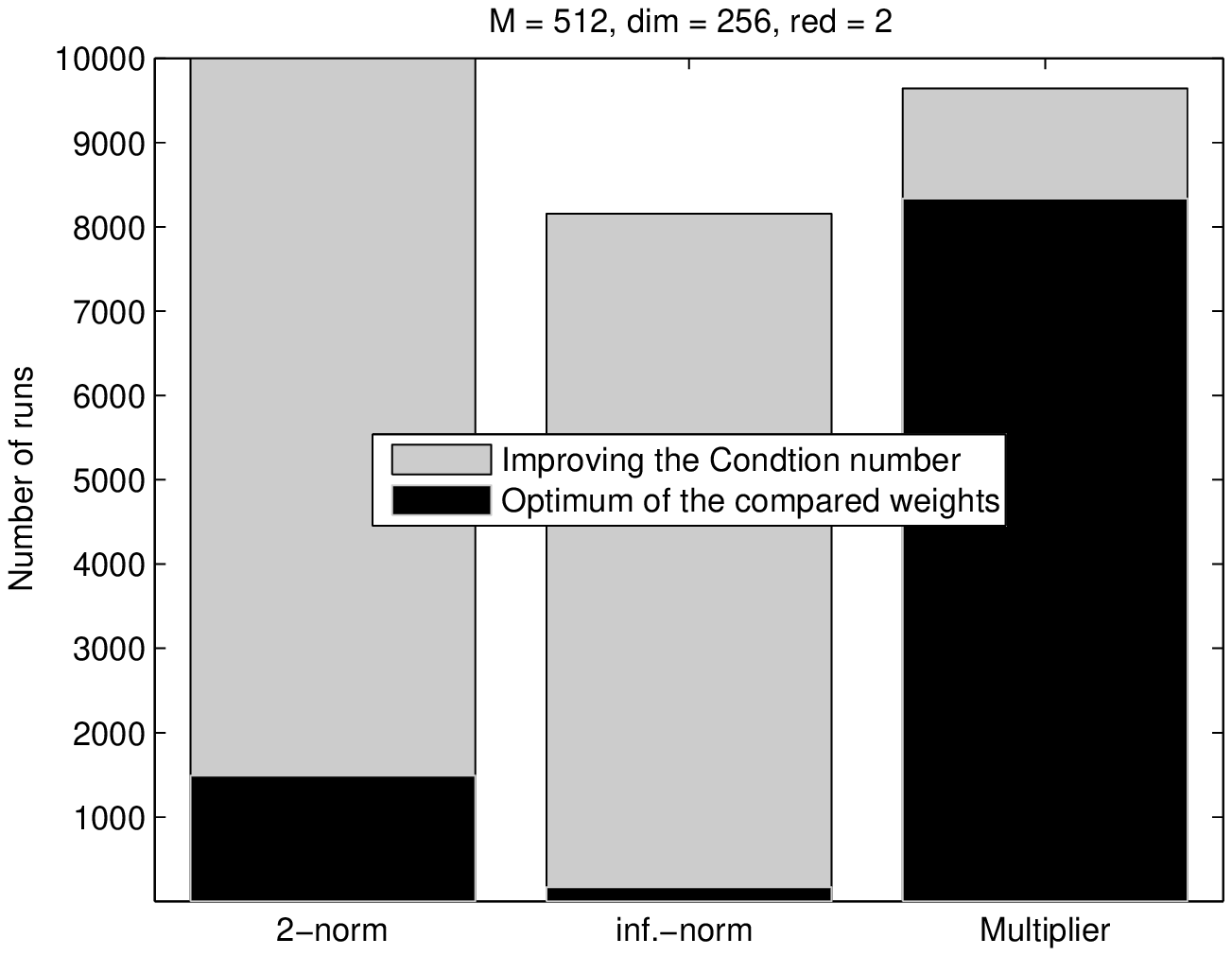}

\includegraphics[width=0.4\textwidth]{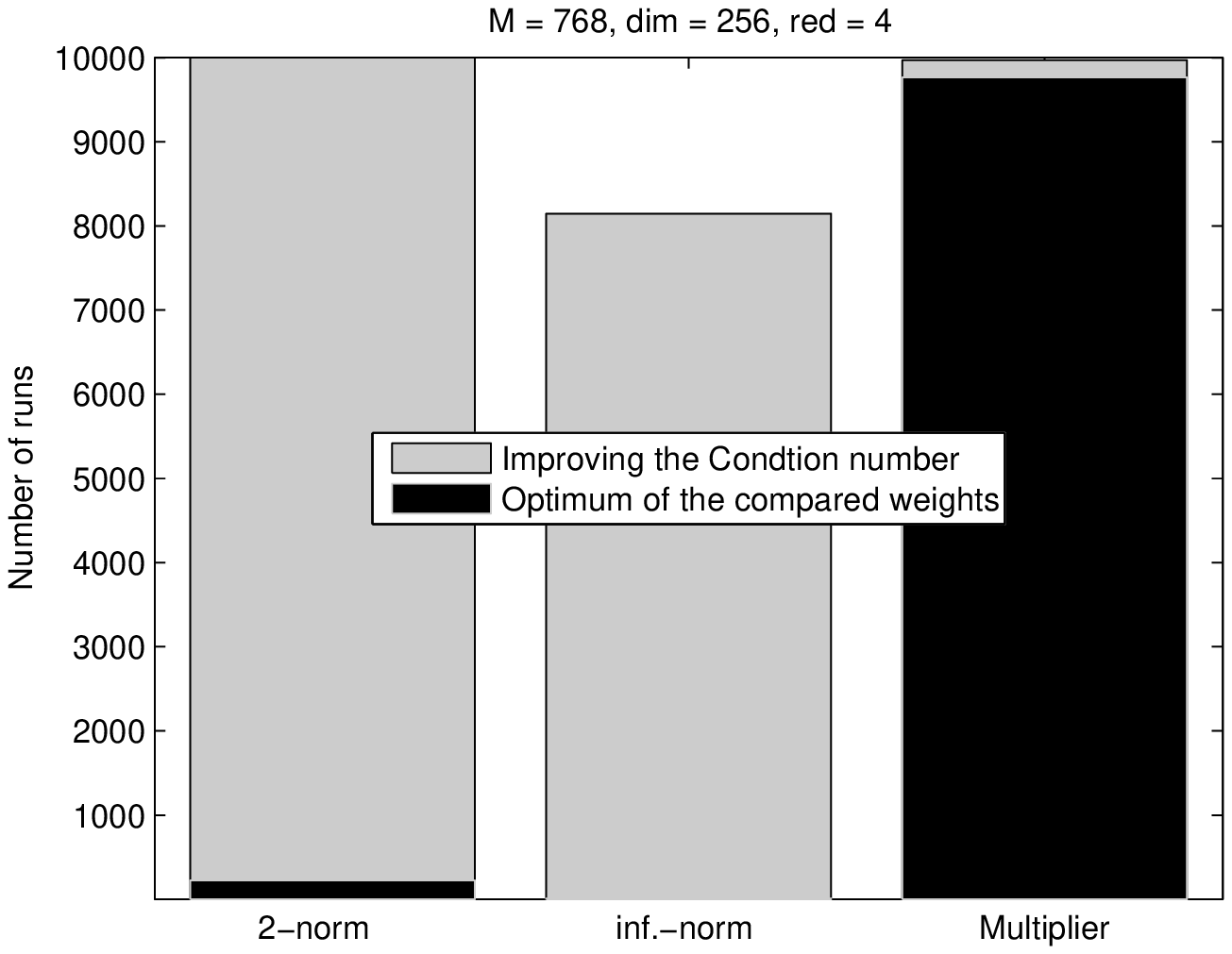}\vspace{3mm}
\includegraphics[width=0.4\textwidth]{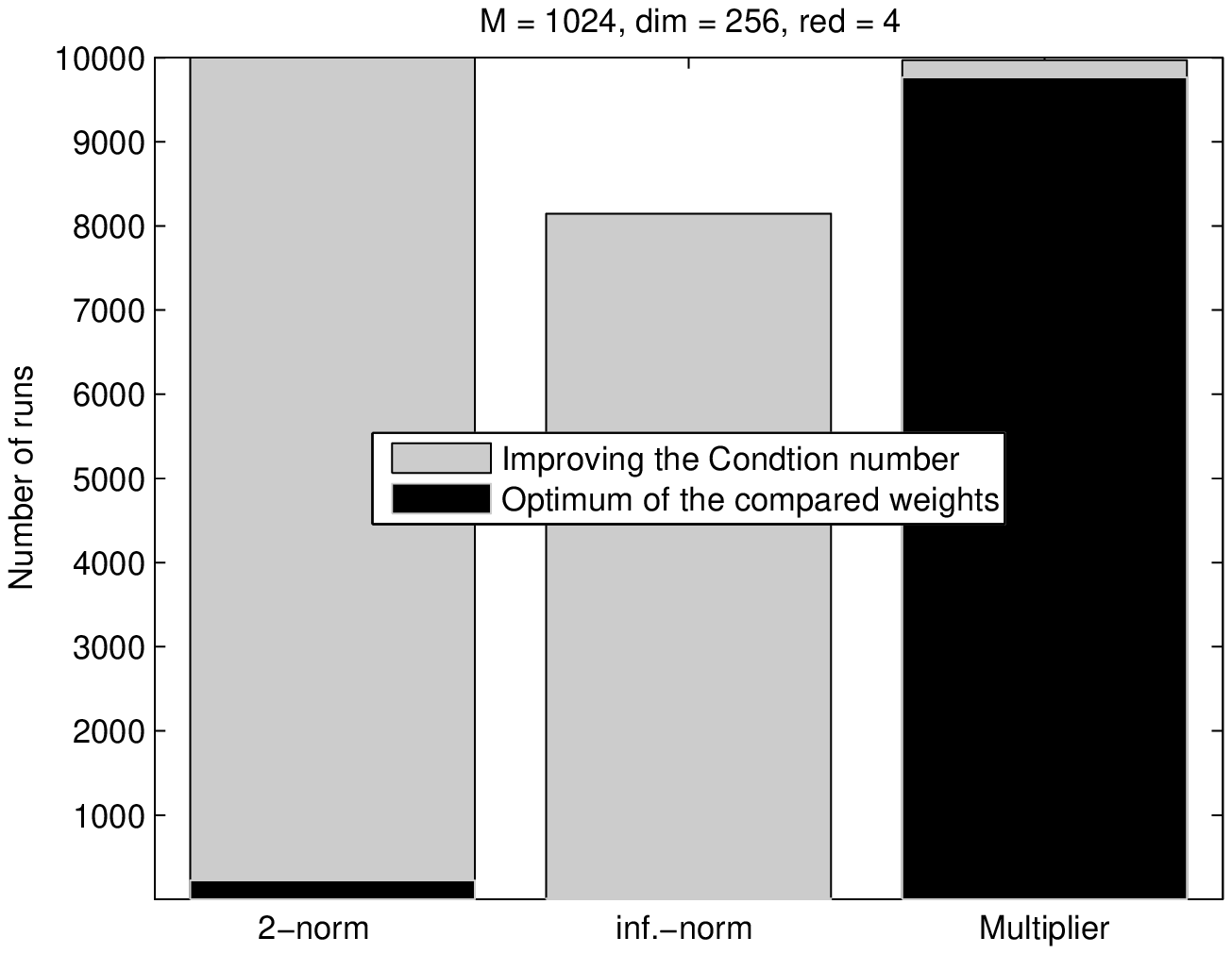}
\end{center}
\caption{\label{fig:testtight4}\small \em  Frames  in $d=256$ dimensions, with the same conventions as in Figure \ref{fig:testtight3}.
Top left: Frame with $M = 260$ elements; top right: $M = 512$; bottom left: $M = 768$; bottom right: $M = 1024$.} 
\end{figure}

In Figures \ref{fig:testtight3} and \ref{fig:testtight4}, we present the results for dimensionality $d = 64$ and $d = 256$, respectively,
and increasing redundancy $(=M/d)$. The graphs show the improvement of condition number by weights $\omega^{(2)}$ (= `2-norm'), 
 $\omega^{(\infty)}$ (= `inf.-norm') and $\omega^{\rm (mult)}$ (= `Multiplier'), respectively. 
 Other tests have been made with lower
dimensionalities $d = 3$ and $d = 10$.  It turns out that, for these cases, the results are sowewhat erratic. Since such low dimensions are not very realistic for applications, 
we simply dropped them. Similarly, 
weights $\omega^{(4)}$ and $\omega^{(6)}$ lead to a higher computational load and give worse results.
Thus these weights are no longer considered.

 Figure \ref{fig:testtight3} shows the results for the parameters: $d = 64$ and $M = 260, 512, 1024$ and $2048$. 
 The condition number is improved in  $96.54 \%$, $99.98\%$, $100 \%$, resp.$100 \%$ of the tests. 

In Figure \ref{fig:testtight4} we have summarized the results for a set of parameters which may be more realistic for  applications, namely, 
$d = 256$ and $M = 260, 512, 1024$ and $2048$. Now the condition number is improved in  $98.87 \%$, $99.97 \%$, $100 \%$, resp. $100 \%$ of the tests. 

%\begin{figure}[ht!]
%\begin{center}
%\includegraphics[width=0.48\textwidth]{condimprov_123_100_new.eps}
%\includegraphics[width=0.48\textwidth]{condimprov_4010_100_new.eps}
%
%\includegraphics[width=0.48\textwidth]{condimprov_25664_100_new.eps}
%\includegraphics[width=0.48\textwidth]{condimprov_1024256_100_new.eps}
%\end{center}
%\caption{\label{fig:testcond} Improvement of condition number by the optimal weight. The original condition numbers are connected by a full line, 
%while the improved condition number are connected by a dotted line. Top Left: $d = 3, M = 12$, top right: $d = 10, M = 40$. Bottom left: $d = 64, M = 256$, 
%bottom right:  $d = 256, M = 1048$} 
%\end{figure}
\vspace*{-2mm}
\subsection{Interpretation} 

The hope, of course, was to find a clear trend with increasing dimensionality and/or redundancy, but  the numerical experiments do not allow such a conclusion.

%First of all, Figures \ref{fig:testtight1} and \ref{fig:testtight2} clearly shows that $\omega^{(p)}$ for $p=4,6$ are almost always inferior to the other weights. 
%Therefore they have not been included in the higher dimensional cases ($d=64$ and $d=256$). In preliminary tests, the same result has been obtained for $p = 1$ and $3$. 

For low dimension ($d=3$ and $d=10$, not shown here) and low redundancy,  the weight $\omega ^{(\infty)}$ and the weight by multiplier approximation $\omega ^{\rm (mult)}$
are sometimes  good, but not always. For high redundancy, $\omega ^{(2)}$ seems to be always a good guess, nearly always improving the condition number. 

%
% does not always improve the condition number, but if it does, 
%it is quite likely to be the optimal of the three presented weights. 
%In dimension $d=10$,  $w^{(\infty)}$ is surpassing the other weights for low redundancy, but for high redundancy $\omega ^{(2)}$ is nearly always the optimal weight. 
%Both in dimensions $d=3$ and $d=10$,  two-fold redundancy  seems to be a special case (for opposite reasons, actually), but this is probably not significant.

For higher dimensions ($d = 64$ and $d=256$), again the weight $\omega ^{(2)}$   nearly always improves the condition number, but, especially for higher redundancy, 
the weight by multiplier is the optimal solution of the three tested weights.

%In order to see how much the condition number is improved in each case, we plot in Figure \ref{fig:testcond} the original condition number 
%and the improved condition number (by the optimal weights). Note that a high redundancy of 4 has been chosen, as in this case the improvement is better. 
%Clearly, there is some improvement in all cases, but the behavior is quite erratic.

As a general rule, however, in order to improve the numerical behavior of frames, the `power weight'  $\omega ^{(2)}$  should be used, because the weight by multiplier approximation is a 
highly complex algorithm and the weighting by the `power weight'   nearly always improves the condition number. Thus $\omega ^{(2)}$ is a good compromise.

However, the only conclusion of these preliminary results is that the connection between optimal weight, dimensionality and redundancy should  be further investigated.

%\bigskip

%%%%%%%%%%%%%%%%%%%%%%%%%%%%%%%%%%%%%%%%%%%%%%%%%%%%%%%%%%%%%%%%%%%%%%%%%%%%%%%%%%%%%%%%%%%%%%%%%%%%%%%%%%%%
\vspace*{-4mm}
\section{Numerical results for discrete Gabor frames} 
\label{sec:gabornum0}

In this last section, we shall examine the case of discrete Gabor frames in concrete situations.
We denote a given Gabor system by $G = \left( \pi ( \lambda) g \right)$ and, for a given weight $(\omega_\lambda)$, the weighted Gabor system by 
\emph{WG} = $\left( \omega_\lambda \pi ( \lambda) g \right)$. Furthermore we will use the notation \emph{DWG} for the canonical dual of the weighted Gabor system, i.e., 
\emph{DWG} = $(\widetilde{\omega_\lambda \pi ( \lambda) g})  %\widetilde{\phantom{w}}
$ and by \emph {iWDG} the dual frame weighted with the reciprocal weights, i.e., 
 \emph {iWDG} = $\big( \frac{1}{\omega_\lambda} \pi ( \lambda) \tilde{g} \big)$.

According to Lemma \ref{sec:framepertmul1}, for semi-normalized weights, \emph {iWDG} is a dual frame of \emph{WG}, but not necessarily the canonical dual. 
In this section we investigate how close these two duals are to each other, i.e.,  how well \emph {iWDG}  approximates \emph{DWG}.
The rationale behind the examples is the following. Among all possible duals, the canonical one is the unique one that satisfies the minimal norm condition.
 However, it is often difficult to compute. On the contrary,  \emph {iWDG} is much easier to evaluate, and thus could be used as a convenient substitute for \emph{DWG}.

We treat the cases with several different windows, in dimension $d=144$. In addition, we consider our frames with $M$ elements as $d$-periodic.
%, that is, each component belonging to $\CC^{d}$%=\{0,1,\ldots,d-1\}$
Explicit results will be given for  a Gaussian, a Hanning and a 
Bartlett window. %(see Fig \ref{fig:threewindBW}). 
Similar computations have been performed also for a Blackman window and B-spline windows of order 3 and 5, 
but the results are not significantly different, so we will skip them here.

%\begin{figure}[!ht]
%\begin{center}
%\includegraphics[width=0.8\textwidth]{threewindBW.eps}
% threewindBW.eps: 300dpi, width=3.79cm, height=2.85cm, bb=   81   227   529   564
%\caption{\label{fig:threewindBW} \small \em The window functions used in the examples are presented in the centered plot. 
%Only the relevant nonzero part is displayed, i.e. supported on [-30,30].}
%          \end{center}
%\end{figure}

\begin{table}[!ht]
 \begin{center} \small
\begin{tabular}{|c|c|c|c|}
\hline
$(a,b)$  & Gaussian & Hanning & Bartlett \\ 
\hline
$(12,9)$ & 2.5041 & 2.8609 & 4.9648  \\
\hline
$(9,8)$  & 1.4258 & 2.0000 & 3.9512   \\
\hline
$(8,6)$  & 1.1324 & 1.1603 & 1.5612  \\
\hline
$(6,6)$  & 1.0151 & 1.1266 & 1.4483  \\
\hline
$(6,4)$  & 1.0075 & 1.0000 & 1.0857  \\
\hline
$(4,4)$  & 1.0000 & 1.0000 & 1.0375  \\
\hline
\end{tabular}
\caption{\label{tab:FrameBounds}\small  Frame bound ratio ({\sf M}/{\sf m}) of the Gabor frame $G$ calculated for the given windows and time-frequency shifts.}
\end{center}
\end{table}
%\vspace*{-2mm}
Our frame elements (atoms) read as $g_{k,l} = M_{lb}T_{ka}g$, where $k=0,1,\ldots, \frac{d}{a} -1$ and $l=0,1,\ldots, \frac{d}{b} -1$. Thus the number of frame elements %dimension of the frame
is $M=rd$, where $r:=d/ab$ is the redundancy.
Six pairs of time-frequency shift parameters are considered  to construct the lattices, namely, (12,9), (9,8), (8,6), (6,6), (6,4), (4,4),  with redundancy 1.33, 2, 3, 4, 6, 9, respectively. 
The frame bound ratio for the Gabor frame $G$ calculated for each of the given window functions and time-frequency parameters are presented in Table \ref{tab:FrameBounds}.

To investigate the error of approximation of \emph{DWG} by \emph{iWDG}, we consider the relative error in Hilbert-Schmidt norm of the 
two synthesis matrices: % namely \texttt{Err=norm(iWDG-DWG)/norm(DWG)}.
$$
 \epsilon = \frac{\norm{\rm HS}{iWDG-DWG}}{\norm{\rm HS}{DWG}} = \sqrt{
\frac{\sum\limits_\lambda \norm{\CC^d}{ \frac{1}{\omega_\lambda} \pi ( \lambda) \tilde{g} - \widetilde{({\omega_\lambda \pi ( \lambda) g})}}^2}
{\sum\limits_\lambda \norm{\CC^d}{\widetilde{({\omega_\lambda \pi ( \lambda) g})}}^2}}.
$$
The formula is related to the notions of `quadratic closeness' \cite{young1}, and  `Bessel norm'.\cite{xxlmult1}

To start with, let us consider a periodized Gaussian window $g$ of length $d=144$ and the lattice parameters $a=12, b=9$, so that $d/a=12, d/b=16$ and $M=192$.  

We consider a weight that
mimicks a local mask, such as the one used in Ref.\cite{jpa-JFAA} for enhancing the contrast in a picture of the Red Spot of Jupiter. 
%Namely, the weight takes the constant value between 1 and 2 everywhere, except for a central block where it equals 2. We also vary the size of that central block.
Namely, we take the weight $w$ equal 2 on the centered $3\times 3$ block on the lattice and equal 1 outside of this block. 
Hence only nine out of the 192 elements of the frame $G$ are amplified, while the rest is unchanged, % (Figure \ref{fig:boxweight})
 i.e., 
$$
 \omega_{k,l} = \left\{ 
\begin{array}{c l} 
2, & (k,l) \in \left\{ \frac{d}{a}-p,\ldots,\frac{d}{a}-1,0,\ldots,p \right\} \times 
\left\{  \frac{d}{b}-p,\ldots,\frac{d}{b}-1,0,\ldots,p \right\} \\
1, & \mathrm{ otherwise } \end{array} \right. 
$$
where $k=0,1,\ldots \frac{d}{a} -1$ and $l=0,1,\ldots \frac{d}{b} -1$ and the parameter $p$ equals 1. 

\begin{figure}[t]
\begin{center}
\includegraphics[width=0.32\textwidth]{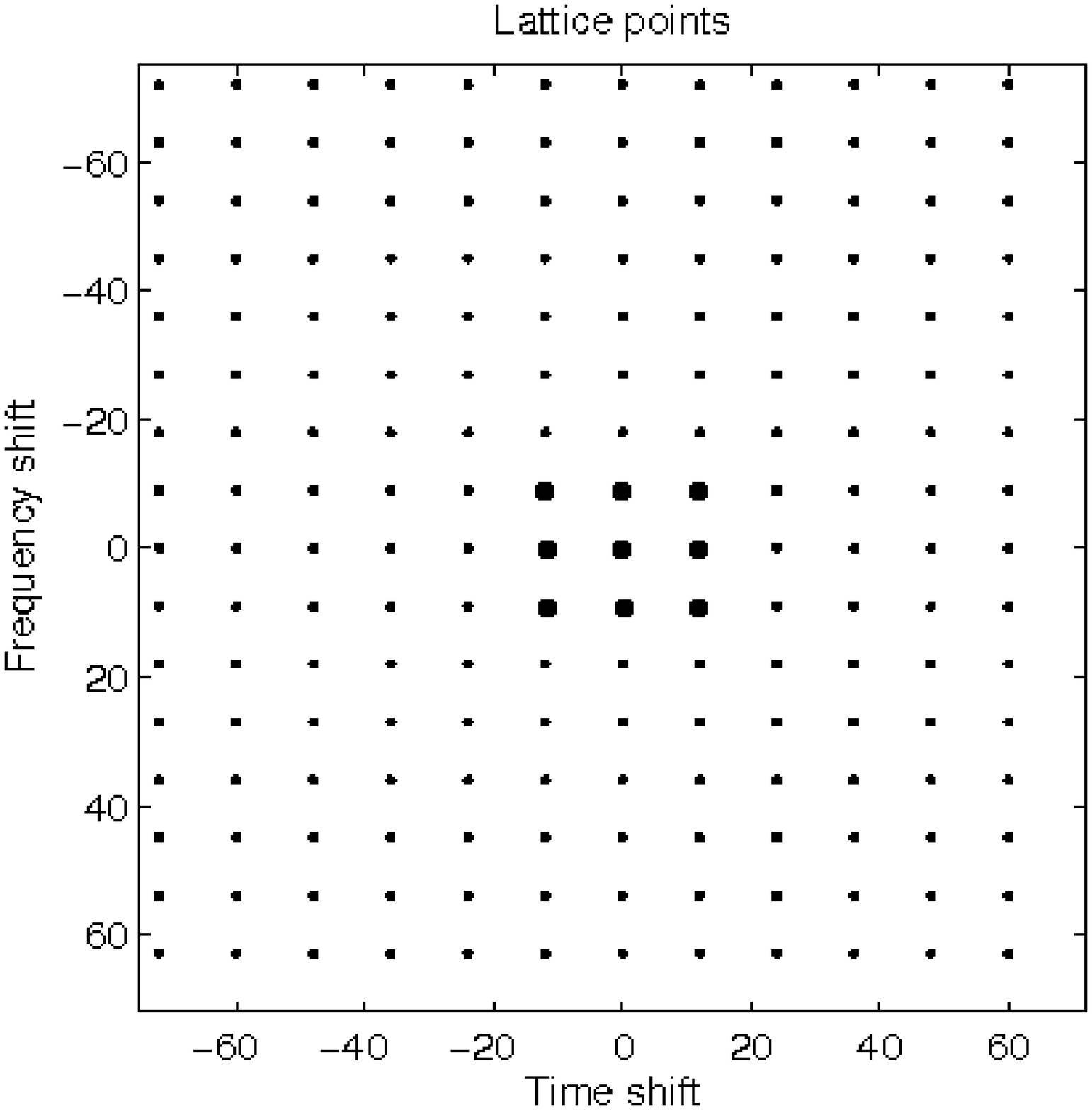}
\includegraphics[width=0.32\textwidth]{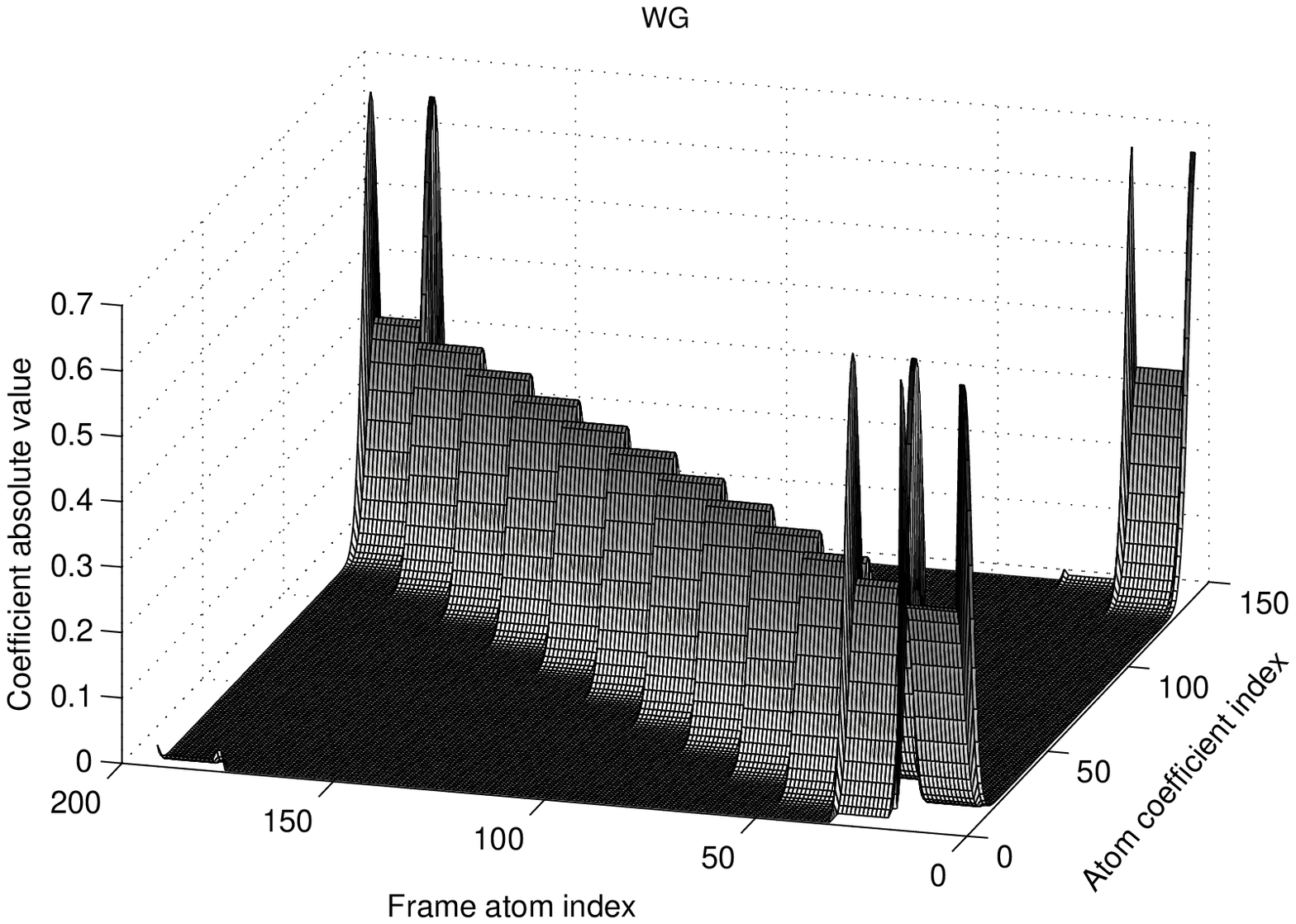}
\includegraphics[width=0.32\textwidth]{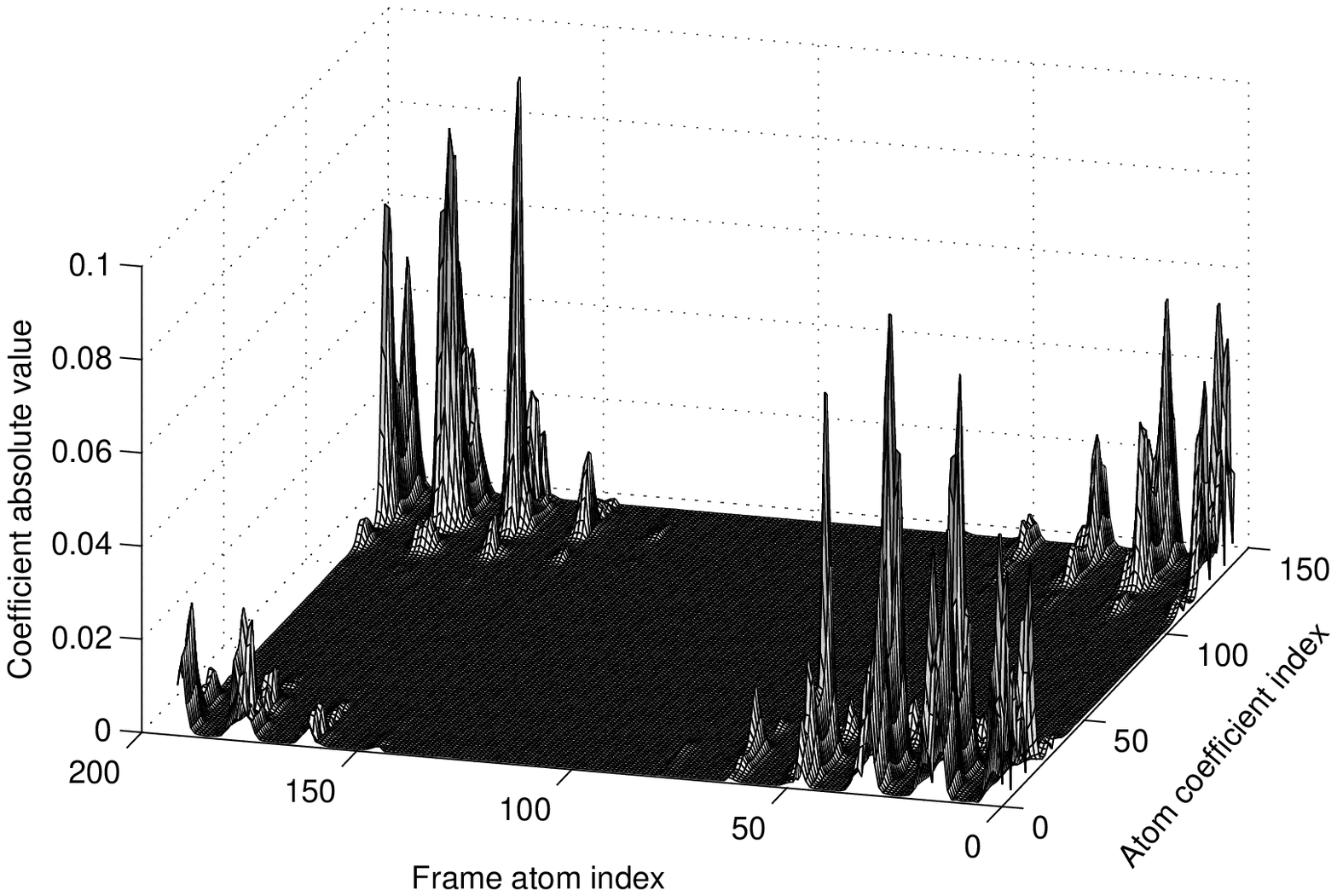}
\caption{\label{fig:boxweight} \small \em Nine elements of the frame $G$ are amplified with the weight equal to 2 while the rest remains unchanged.
 (Left) Positions of the amplified atoms in time-frequency domain (the plot is rescaled in order to group the marked points in the center). 
(Center) The resulting weighted frame WG.  Note the ``spikes" which are the amplified atoms around (0,0).
%}
%\end{center}
%\end{figure}
%\begin{figure}[!ht]
%\begin{center}
%\caption{\label{fig:diffbox} \small \em 
(Right) Difference between \emph{iWDG} and  \emph{DWG}.
The most notable changes are located near the places where the weight $\omega=2$ was applied.
}
\end{center}
\end{figure}

The weighted frame $WG$ is presented in Figure \ref{fig:boxweight} in such a way that the atoms $g_{k,l}$ (each with 144 coefficients) are stacked 
along the $y$ axis (``frame atom index"). The actual ordering of the frame atoms is arbitrary, because it does not influence the frame operator $L$. 
Here atoms with the same  time shift, but different modulations, are stacked one alongside the other, starting with the smallest time shifts. 

Since the weights used are positive numbers, it suffices to display the absolute values of the atom coefficients ($z$ axis on Figure \ref{fig:boxweight})
 in order to show the effect of the weights. 
However, as a side effect, the atom modulation is not visible anymore. Atoms which differ only in modulation are grouped in clearly visible bands 
containing $d/b=16$ atoms each.

%\begin{figure}[!ht]
%\begin{center}
%\includegraphics[width=0.7\textwidth]{sinweightGray.eps}
%\caption{\label{fig:sinweight} \small \em Gabor frame weighted with a sinusoidal weight: Absolute values of the frame elements.}. 
%          \end{center}
%\end{figure}

%\begin{figure}[!ht]
%\begin{center}
%\includegraphics[width=0.7\textwidth]{dualsinGray.eps}
%\caption{\label{fig:dualsin} \small \em Canonical dual frame of the (sinusoidal) weighted Gabor frame from Figure \ref{fig:sinweight}: 
%Absolute values of the frame elements. }
% \end{center}
%\end{figure}

%\begin{table}[!ht]
%  \begin{center}
%\begin{tabular}{|c|c|c|c|c|c|c|}
%\hline
%$(a,b)$  & Gaussian & Hanning & Bartlett  \\ 
%\hline
%$(12,9)$ & 0.1829   & 0.2775   & 0.3441    \\
%\hline                         
%$(9,8)$  & 0.1601   & 0.2550   & 0.3121    \\
%\hline                         
%$(8,6)$  & 0.1574   &  0.2360   & 0.2680   \\
%\hline                         
%$(6,6)$  & 0.1566   &  0.2360   & 0.2676   \\
%\hline                         
%$(6,4)$  & 0.1566   & 0.2353   & 0.2630    \\
%\hline                         
%$(4,4)$  & 0.1566   & 0.2353   & 0.2628    \\
%\hline
%\end{tabular}
%\end{center}
%  \caption{\label{tab:errappdwg1}\small  The relative error of approximation of \emph{DWG} by \emph {iWDG} with respect to the
%Hilbert-Schmidt norm, in the case of Gabor frames with increasing redundancy weighted with a sinusoidal weight.}
%\end{table}

%%%%%%%%%%%%%%%%%%%%%%%%%%%%%%%%%%%%%%%%%%%%%%%%%%%%%%%%

%\begin{figure}[!ht]
%\includegraphics[width=0.45\textwidth]{invdualboxGray.eps}
%\includegraphics[width=0.45\textwidth]{candualboxGray.eps}
%\caption{\label{fig:boxduals} \small \em Comparison between the inversely weighted canonical dual frame \emph{iWDG} (on the left) 
%and canonical dual frame \emph{DWG} 
%to the weighted frame \emph{WG} on the right. 
%.}
%\end{figure}

It turns out that the inversely weighted canonical dual frame 
(\emph{iWDG}) and the canonical dual weighted frame (\emph{DWG}) yield almost  identical plots. 
Thus we present in Figure \ref{fig:boxweight} the absolute value of the componentwise difference \emph{iWDG-DWG} between the two dual families. 
The only visible differences are in the locations where the weight $\omega =2$ was applied (see Figure \ref{fig:boxweight}).

Now we construct frames with higher redundancy by using the time-frequency shift parameters listed above and the different window functions.
The same weight $w$ is applied to the new frames. 
The relative error of approximation of \emph{DWG} by \emph {iWDG}, measured in the Hilbert-Schmidt norm, is presented in Table \ref{tab:errtabbox}.

\begin{table}[!ht]
  \begin{center} \small
\begin{tabular}{|c|c|c|c|c|c|c|}
\hline
$(a,b)$  & Gaussian & Hanning & Bartlett \\ 
\hline
$(12,9)$ & 0.0802 &  0.0830 & 0.0909   \\
\hline                               
$(9,8)$  & 0.0808 &  0.0858 & 0.0908   \\
\hline                               
$(8,6)$  & 0.0786 &  0.0798 & 0.0816   \\
\hline                               
$(6,6)$  & 0.0751 &  0.0779 & 0.0796   \\
\hline                               
$(6,4)$  & 0.0703 &  0.0707 & 0.0718   \\
\hline                               
$(4,4)$  & 0.0638 &  0.0665 & 0.0680   \\
\hline
\end{tabular}
\caption{\label{tab:errtabbox}\small  The relative error of approximation of \emph{DWG} by \emph{iWDG} with respect to the Hilbert-Schmidt norm, 
in the case of Gabor frames weighted with the piecewise constant weights and a $3\times 3$ block on the lattice.
 The redundancy  increases with decreasing $a$ and $b$.}
\end{center}
\end{table}

\vspace*{-2mm}
Next, let us change the number of amplified elements by increasing the size 
of the central square block on the lattice. 
The weight $w$  equals  2 on this block and  1 outside. The size of the block (mask)  changes from $3\times 3, 5\times 5$
to $7 \times 7$ and $9 \times 9$, hence the number of  amplified elements of the frame $G$ increases from 9, 25, to 49 and 81.
For the   case of a Gabor frame with a Gaussian window $g$ and time-frequency parameters $(a,b)=(12,9)$, this constitutes 4.7\%, 13\%, 25.5\% 
and 42.2\% of all frame elements, respectively. We apply the weights $ \omega_{k,l}$ with the parameter $p=2$ for the $5\times 5$ block, $p=3$ 
for the $7\times 7$ block and $p=4$ for the $9 \times 9$ block.

\begin{figure}[!ht]
\begin{center}
\includegraphics[width=0.4\textwidth]{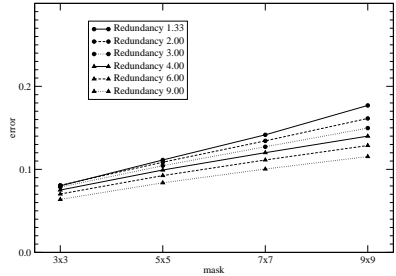}
\caption{\label{fig:changeblock} \small \em Linear dependence of the relative error of approximation on the mask size for varying redundancies. 
Increased mask size leads to a 
larger number of amplified frame elements and  generates a larger error when approximating the canonical dual frame  \emph{DWG}
with the inversely weighted dual frame \emph{iWDG}. Higher redundancy leads to lower errors and this effect is stronger for larger masks.
}
\end{center}
\end{figure}

Figure \ref{fig:changeblock} shows the results for the Gaussian window, for the other window functions the results are very similar, hence we skip them here.
The relative error of approximation of \emph{DWG}  by \emph{iWDG} with respect to the Hilbert-Schmidt norm is presented as a function of the increasing 
block size. The simulations are repeated for the frames with higher redundancy. In each case, the error increases linearly with the size of the block on the lattice. 
The linear relation is observed also for the other window functions.

\vspace*{-4mm}
\subsection{Interpretation}

The redundancy of the frame has a stronger influence on the error of the approximation of \emph{DWG} by  \emph{iWDG} than the size of the applied mask. 
Interestingly the influence is almost linear.

Another important factor is the choice of the window, but Gabor frames exhibit a remarkable indifference to the change of window in the situations we considered. 
Although the Gaussian window gave the best results, other window functions such as Hanning, Bartlett, Blackmann or B-spline windows of order 3 and 5 
are only slightly worse.
\vspace*{-4mm}
\section{Perspectives}

For finding the optimal weights, further numerical tests should be conducted. In particular, a geometric classification of those frames, where the given
  weights work well or not at all, will be interesting. Furthermore, the Hilbert-Schmidt operator norm is easy to use for measuring the approximation error,
 but it would be more interesting to use  the operator norm. This could be obtained, for instance, by using LMI algorithms.\cite{boyd1} 

The weights in Section \ref{sec:weightnum0} can be seen as a measure of the off-diagonal behavior of the Gram matrix. It is well-known\cite{forngroech1} that the 
cross-Gram matrix 
of the frame and its canonical dual, i.e., $G_{\Psi, \widetilde\Psi}$, carries a lot of information about the frame. 
In further numerical tests, we will evaluate the weights using this matrix.
%
%Finally, it is known \cite{xxlmult1} that, in the general frame case, no exact symbolic calculus for multipliers can be assumed, 
%i.e., the combination of symbols does not always correspond to the combination of the operators. 
%So even in the case of a single frame $(\psi_n)$, we get in general
%\begin{equation} \label{sec:combmult1} 
%{\bf M}_{\bf m^{(1)},\Psi} \circ {\bf M}_{\bf m^{(2)},\Psi} \not= {\bf M}_{\bf m^{(1)} \cdot m^{(2)} , \Psi} .
%\end{equation}
%In general the product of two frame multipliers is not a frame multiplier any more. But it is also known that, if biorthogonal frames are used 
%for analysis and synthesis,  then 
% \eqref{sec:combmult1} becomes an equality.
%
%It is easy to show that the set of multipliers forms a vector space, but with the above properties we know that it normally does \emph{not} form an algebra.
% So the space of the frame multipliers is an example of a partial algebra \cite{antino1}, of a type different from all the known ones. 
%As such it deserves a systematic study,  
%that will be the subject of future work.
%\vspace*{-5mm}
\section*{Acknowledgement}
The authors would like to thank K. Schnass for helpful comments and suggestions.  
This work was partly supported by 
the %European Union's Human Potential Program, under contract HPRN-CT-2002-00285 (HASSIP).
WWTF project MULAC (`Frame Multipliers: Theory and Application in Acoustics', MA07-025).


{  \small
\begin{thebibliography}{10}

%\bibitem{antino1}
%J-P. Antoine, A.~Inoue, and C.~Trapani,
%\newblock {\em Partial *-Algebras and Their Operator Realizations},  
%  Mathematics and Its Applications, vol. 553,
%\newblock Kluwer, Dordrecht, NL, 2002.

\bibitem{jpa-JFAA}J-P.~Antoine and P.~Vandergheynst, Wavelets on the two-sphere and other conic sections,
\emph{J. of Fourier Analysis and Appl.} {13}:   {369--386} {(2007)}.

\bibitem{xxlmult1}
P.~Balazs,
\newblock Basic definition and properties of {B}essel multipliers,
\newblock {\em J. Math. Anal. Appl.},  325(1): 571--585  (2007).

\bibitem{xxlfinfram1}
P.~Balazs.
\newblock Frames and finite dimensionality: Frame transformation,  classification and algorithms.
\newblock {\em Applied Mathematical Sciences}, 2(41--44):2131--2144 (2008).

\bibitem{xxlframehs07}
P.~Balazs,
\newblock {H}ilbert-{S}chmidt operators and frames - {C}lassification, best
  approximation by multipliers and algorithms,
\newblock {\em Int. J. Wavelets, Multires. and Inform. Proc}, 6(2): 315 -- 330  (2008).

\bibitem{xxlframoper1}
P.~Balazs,
\newblock Matrix-representation of operators using frames,
\newblock {\em Sampling Th. Sign. Image Proc. (STSIP)},
  7(1): 39--54 (2008).

\bibitem{xxlfei1}
P.~Balazs, H.~G. Feichtinger, M.~Hampejs, and G.~Kracher,
\newblock Double preconditioning for {G}abor frames,
\newblock {\em IEEE Trans.  Signal Process.}, 54(12): 4597--4610 (2006).
 

\bibitem{bogdvan1}
I.~Bogdanova, P.~Vandergheynst, J-P. Antoine, L.~Jacques, and M.~Morvidone,
\newblock Stereographic wavelet frames on the sphere,
\newblock {\em Applied Comput. Harmon. Anal.} , 19: 223--252 (2005).

\bibitem{boyd1}
S.~Boyd, L.~El Ghaoui, E.~Feron, and V.~Balakrishnan,
\newblock {\em Linear Matrix Inequalities in System and Control Theory}.
\newblock SIAM, Philadelphia, PA, 1994.

\bibitem{Casaz1}
P.G.~Casazza,
\newblock The art of frame theory,
\newblock {\em Taiwanese J. Math.}, 4(2): 129--202 (2000).

\bibitem{ole1}
O.~Christensen,
\newblock {\em An Introduction to Frames and {R}iesz Bases},
\newblock Birkh{\"a}user, Basel, Boston, Berlin, 2003.

\bibitem{conw1}
J.B.~Conway,
\newblock {\em A Course in Functional Analysis}.
\newblock  Springer, New York, 2nd edition, 1990.

\bibitem{daubech1}
I.~Daubechies.
\newblock {\em Ten Lectures On Wavelets}.
\newblock CBMS-NSF Regional Conference Series in Applied Mathematics. SIAM
  Philadelphia, 1992.

\bibitem{doerf1}
M.~D{\"o}rfler,
\newblock {\em Gabor Analysis for a Class of Signals Called Music}.
\newblock PhD thesis, University of Vienna, 2003.

\bibitem{forngroech1}
M.~Fornasier and K.~Gr{\"o}chenig.
\newblock Intrinsic localization of frames.
\newblock {\em Constructive Approximation}, 22(3):395--415 (2005).

\bibitem{gohbgol1}
I.~Gohberg, S.~Goldberg, and M.A.~Kaashoek,
\newblock {\em Classes of Linear Operators, vol. I},
\newblock Birkh{\"a}user, Basel,Boston, Berlin, 1990.

\bibitem{Groech1}
K.~Gr{\"o}chenig.
\newblock {\em Foundations of Time-Frequency Analysis}.
\newblock Birkh{\"a}user Boston, 2001.

\bibitem{luen}
D.~Luenberger,
\newblock {\em Linear and Nonlinear Programming},
\newblock Addison-Wesley, Reading, MA, 1984.

\bibitem{pewa02}
I.~Peng and S.~Waldon,
\newblock Signed frames and {H}adamard products of {G}ram matrices,
\newblock {\em Linear Algebra Appl.}, 347: 131--157 (2002).

\bibitem{trebau1}
L.N.~Trefethen and D.~Bau III,
\newblock {\em Numerical Linear Algebra},
\newblock SIAM, Philadelphia, PA, 1997.

\bibitem{Waldron03}
S.~Waldron,
\newblock Generalized {W}elch bound equality sequences are tight frames, 
\newblock {\em IEEE Transactions on Information Theory}, 49(9): 2307--2309 (2003).

\bibitem{young1}
R.M.~Young.
\newblock {\em An Introduction to Nonharmonic Fourier Series}.
\newblock Academic Press, London, 1980.

\end{thebibliography}
}
\end{document}